\input amssymb.sty
\documentclass{amsproc}
\usepackage{graphicx}
\newtheorem{theorem}{Theorem}[section]
\newtheorem{lemma}[theorem]{Lemma}
\newtheorem{prop}[theorem]{Proposition}
\newtheorem{cor}[theorem]{Corollary}

\newtheorem{definition}[theorem]{Definition}
\theoremstyle{definition}

\newcommand{\pa}{\mathcal P_{\alpha}}
\newcommand{\pb}{\mathcal P_{\beta}}
\newcommand{\pw}{\mathcal P_{w}}
\newcommand{\da}{\mathcal D_{\alpha}}
\newcommand{\db}{\mathcal D_{\beta}}

\newcommand{\T}{\mathcal T}
\newcommand{\PP}{\mathcal P}
\newcommand{\qa}{\mathcal Q_{\alpha}}

\newcommand{\q}{\mathcal Q}
\newcommand{\D}{\mathcal D}

\newcommand{\ewo}{\mathcal E_{w_0}}
\newcommand{\ew}{\mathcal E_{w}}
\newcommand{\ea}{\mathcal E_{\alpha}}
\newcommand{\pwo}{\mathcal P_{w_0}}

\theoremstyle{remark}
\newtheorem{remark}[theorem]{Remark}

\numberwithin{equation}{section}

 \DeclareMathSymbol{\vartriangle}{\mathord}{AMSa}{"4D}
\DeclareMathSymbol{\triangledown}{\mathord}{AMSa}{"4F}
\def\tr{\,{\vartriangle}\,}
\def\ut{\,{\triangledown}\,}

\begin{document}

\title{Littelmann paths and Brownian paths}

\author{Philippe Biane}
\address{CNRS, D\'epartement de Math\'ematiques et Applications,
    \'Ecole Normale Sup\'erieure, 45, rue d'Ulm 75005 Paris, FRANCE
}
\email{Philippe.Biane@ens.fr}
\author{ Philippe Bougerol}
\address{Laboratoire de Probabilit\'es et mod\`eles al\'eatoires,  
Universit\'e Pierre et Marie Curie, 4, Place Jussieu, 75005 Paris,  
FRANCE}
\email{bougerol@ccr.jussieu.fr}
\author{ Neil O'Connell}
\address{Mathematics Institute
University of Warwick
Coventry CV4 7AL UK}
\email{noc@maths.warwick.ac.uk}
\subjclass{Primary ; Secondary }
\date{}

\begin{abstract}
We study some path transformations related to Pitman's theorem on Brownian
motion and the three dimensional Bessel process. We relate these to
Littelmann path model, and give 
applications to representation theory and to Brownian motion in a Weyl
chamber.
\end{abstract}
\maketitle
\section{Introduction}
Some transformations defined on continuous paths with values in a
vector
space have appeared in recent years, in two separate parts of
mathematics. On
the one hand Littelmann \cite{littel} developed his path model in
order to give a unified
combinatorial setup for representation theory, generalizing the theory
of Young
tableaux to semi-simple or Kac-Moody Lie algebras of type other than
$A$. On the
other hand,  in probability theory, several path transformations have
been
introduced that yield a construction of Brownian motion in a Weyl
chamber
starting from a Brownian motion in the corresponding Cartan Lie algebra.
The oldest and simplest of these transformations comes from Pitman's
theorem
    \cite{pitman}  which states that if $(B_t)_{t\geq 0}$ is a
one-dimensional Brownian motion,
then the stochastic process $R_t:=B_t-2\inf_{0\leq s\leq t}B_s$ is a
three
dimension Bessel process, i.e. is distributed as the euclidean norm of
a three
dimensional Brownian motion (actually Pitman stated his theorem with the
transformation $2\sup_{0\leq s\leq t}B_s-B_t$, but thanks to the symmetry of
Brownian motion this is clearly equivalent to the above statement). 
It turns out that  the fact that,
here, the dimension of the Brownian motion is equal to  1,  the
rank of the group $SU(2)$, while 3, the dimension of the Bessel
process, is the
dimension of the group $SU(2)$ is not a mere coincidence but a
fundamental fact
which we will clarify in the following. Pitman's theorem has been
extended in
several ways. The first step has been the result of Gravner, Tracy and
Widom,
\cite{GTW} and of Baryshnikov \cite{Bary} which states that the largest
eigenvalue of a random $n\times n$ Hermitian matrix in the GUE
    is distributed as the random variable
$$\sup_{1=t_n\geq t_{n-1}\geq\ldots\geq t_1\geq t_0=0
}\ \sum_{i=1}^n (B_i(t_i)-B_{i}(t_{i-1}))$$
where $(B_1,\ldots,B_n)$ is a standard $n$-dimensional Brownian motion.
This result in turn was  generalized in \cite{bj} and \cite{oy2}. 
These extensions involve path transformations
which generalize Pitman's  and
are closely related to the Littelmann path model.
    One of the purposes of this
paper is to clarify these connections as well as to settle a number of
questions
raised in these works. In the course of these investigations we will
derive 
several applications to representation theory.
    These path transformations occur in
    quite different contexts, since the one in \cite{bj} is expressed by
    representation theoretic means, whereas the one in \cite{oy2} 
    is purely combinatorial, and arises
     from queuing theory considerations.

    Let us describe more precisely the content of the paper. We start by
defining
    the Pitman transforms which are the main object of study in this
paper.
    These transforms operate on the set of continuous functions  
$\pi:[0,T] \to V$, with
values in some
    real vector space $V$, such that $\pi(0)=0$. They
     are given by the formula
    $$\pa \pi(t)=\pi(t)-\inf_{t\geq s\geq 0}\alpha^{\vee}(
\pi(s))\alpha,\qquad
 t\in [0,T].$$
Here $\alpha\in V$ and $\alpha^{\vee}\in V^{\vee}$ (where $V^{\vee}$ is
the dual space
of $V$)  satisfy $\alpha^{\vee}(\alpha)=2$. These are multidimensional
generalizations of the transform occuring in Pitman's theorem.
They are related to Littelmann's operators as shown in section
\ref{littelpit}.
    We show that these
transforms satisfy braid relations, i.e.
if $\alpha,\beta\in V$ and   $\alpha^{\vee},\beta^{\vee}\in V^{\vee}$
     are  such that
$\alpha^{\vee}(\alpha)=\beta^{\vee}(\beta)=2$, and
$\alpha^{\vee}(\beta)< 0,
\beta^{\vee}(\alpha)< 0$ and
$\alpha^{\vee}(\beta)
\beta^{\vee}(\alpha)=4\cos^2\frac{\pi}{n}$, where $n\geq 2$
    is some  integer,
then one has
    $$\pa\pb\pa \ldots=\pb\pa\pb\ldots$$ where there are $n$
factors in each product.
    Consider now a
Coxeter system $(W,S)$ (cf \cite{bo},\cite{humphreys}). 
To each fundamental reflection $s_i$ we associate
a Pitman transform $\mathcal P_{\alpha_i}$. The braid
    relations imply that if $w\in W$ has a reduced
    decomposition $w=s_{1}\ldots s_{n}$,
    then the operator
    $\mathcal P_w=\mathcal P_{\alpha_{1}}\ldots \mathcal
P_{\alpha_{n}}$ is
    well defined, i.e. it depends only on $w$ and not on the reduced
decomposition.
    We show that if $W$ is a Weyl group, $w_0\in W$ is the longest
element, and
    $\pi$ is a dominant path ending in the weight lattice, then for any
path
    $\eta$ in the Littelmann module generated by $\pi$, one has
    \begin{equation}\label{Bpi}
    \pi=\mathcal P_{w_0}\eta.
    \end{equation}
    The path transformation introduced in \cite{oy2}
 can be expressed as $\mathcal P_{w_0}$ where $w_0$ is the
longest
    element in the Coxeter group of type $A$.

    We derive a representation theoretic formula for $\mathcal P_{w}$, in
    the case of a Weyl group,
    expressed in terms of
    representations of the  Langlands dual group, see Theorem \ref{bjoc}.
    This formula is canonical, in the sense that it is
    independent of any choice of a reduced decomposition of
    $w$ in the Weyl group. It is obtained by lifting the path to a path
$g(t)$ with
    values
    in the
    Borel subgroup of the simply connected complex Lie group associated
with the
    root system. Then one obtains integral transformations which relate
the
    diagonal parts in the Gauss decompositions of the elements $\overline
wg(t)$. The
    Pitman transforms are obtained by going down to the Cartan algebra by
applying
    Laplace's method.
    By (\ref{Bpi}) we obtain in this way  a new formula for the
    dominant path in some Littelmann module, in terms of any of the paths
of the
    module, which is a  generalization to arbitrary root systems 
     of Greene's formula (see \cite{ful}). As a byproduct of this formula we also
    obtain a direct proof of the symmetry of the Littlewood-Richardson
    coefficients.

    This formula appeared in \cite{bj} where it was conjectured that the
associated
    map transforms
    a Brownian motion in the Cartan Lie algebra into a Brownian motion in
the Weyl
    chamber. This conjecture was proved in \cite{bj} for some classical  
groups.
    Here we give a completely different proof, valid for all root  
systems.

    This paper is organized as follows. In section 2 we define the
    elementary Pitman transformations
     operating on continuous paths with values in some real
    vector
    space $V$, taking the value $0$ at $0$. The first result is a formula
for the
repeated compositions of  two Pitman
transforms, which implies that they satisfy the
    braid relations. Then we  define  Pitman transformations
     $\PP_w$ associated
to a Coxeter  system $(W,S)$.
In section 3 we prove our main result
which is a
representation theoretic formula for these operators $\PP_w$ in the
case where $W$
is a Weyl group. This formula unifies the results of \cite{oy2}  and of \cite{bj}.
 Results  of  Berenstein  
and Zelevinsky \cite{beze}
and of
Fomin and Zelevinsky \cite{foze}
on totally positive matrices play a crucial role in the proof.
In section 4 we make some comments on a duality transformation  
naturally
defined on paths, which generalizes the Sch\"utzenberger involution, and give an
application to the symmetry of the Littlewood-Richardson rule.
In section 5 we give two proofs of
the generalization of the representation of
Brownian motion in a Weyl chamber obtained in \cite{oy2} and \cite{bj}.
One of the proofs relies essentially
on the duality
properties,  while the other uses Littelmann paths in the context of  
Weyl groups.
Finally section 6 is an appendix where we have postponed
a technical
proof.
\medskip

{\it Acknowledgements. }
We would like to thank P. Littelmann for a 
 useful conversation at an early stage of this work, 
and P. Diaconis and S. Evans for helpful discussions.
   We also thank the referee for useful comments.

\section{Braid relations for the Pitman transforms}
\subsection{Pitman transforms}
Let $V$ be a real vector space, with dual space $V^{\vee}$.
Let $\alpha\in V$
and $\alpha^{\vee}\in V^{\vee}$ be such that $\alpha^{\vee}(\alpha)=2$.
\begin{definition}\label{pitman-transform}The Pitman
    transform  $\pa$
    is defined on the set of continuous paths $\pi:[0,T]\to V$,  
satisfying
     $\pi(0)=0$, by the formula:
$$\pa \pi(t)=\pi(t)-\inf_{t\geq s\geq 0}\alpha^{\vee}(
\pi(s))\alpha,\qquad
T\geq t\geq 0.$$
\end{definition}
This transformation seems to have appeared for the first time in
\cite{pitman}
in the one-dimensional case.
Note that $\pa$ actually depends on the pair $(\alpha,\alpha^{\vee})$.
For
simplicity we shall use the notation $\pa$, it will be always clear
from the
context which $\alpha^{\vee}$ is involved.

When, for some $v \in V$, $\pi$ is the linear path $\pi(t)=tv$ then  
$\PP_\alpha\pi=\pi$ when $\alpha^{\vee}(
v) \geq 0$ and $\PP_\alpha\pi=s_\alpha\pi$ when $\alpha^{\vee}(
v) \leq 0$ where $s_\alpha$ is the reflection on $V$
\begin{equation}\label{salpha}s_{\alpha} v=v-\alpha^{\vee}(
v)\alpha
\end{equation}
for $v \in V$.

We list a number of elementary properties of the Pitman transform below.
\begin{prop}\label{pit}

({\it i}) For any  $\lambda>0$  the Pitman transformation associated
with
the pair $(\lambda\alpha,\alpha^{\vee}/\lambda)$
is the same as the one associated with the pair
$(\alpha,\alpha^{\vee})$.

({\it ii}) One has $\alpha^{\vee}(\pa \pi(t))\geq 0$ for all $t\in
[0,T]$.
Furthermore
    $\pa \pi=\pi$ if and only if  $\alpha^{\vee}(\pi(t))\geq 0$
for all
$t\in [0,T]$.

({\it iii} )
The transformation $\pa$ is an idempotent, i.e. $\pa\pa \pi=\pa \pi$
for all
$\pi$.

({\it iv})
Let $\pi:[0,\infty[\to V$ be a path, then
     $-\inf_{0\leq t\leq T}\alpha^{\vee}(\pi(t))\in [0,\alpha^{\vee}(
\pa\pi(T) )]$.
     Conversely, given a path $\eta$
      satisfying $\eta(0)=0$,
    $\alpha^{\vee} (\eta(t))\geq 0$
for all
$t\in [0,T]$ and $x\in [0,\alpha^{\vee}( \eta(T) )]$,
     there exists a unique path $\pi$ such that
$\pa \pi=\eta$ and
    $x=-\inf_{T\geq t\geq 0}\alpha^{\vee}(\pi(t))$.
    Actually $\pi$ is given by the
formula
\begin{equation}\label{eta-pi}
\pi(t)=\eta(t)-\min\left(x,\inf_{T\geq s\geq t}
\alpha^{\vee}( \eta(s) )\right)\alpha.
\end{equation}

\end{prop}
{ Proof.}
Items ({\it i}) and ({\it ii}) are trivial, and ({\it iii}) follows
immediately
from ({\it ii}).
    Hopefully the reader can
give a formal proof of ({\it iv}), see section \ref{prfpit} for such a
proof, but
it is perhaps more illuminating to
stare for a few minutes
at Fig. 1, which shows, in the one dimensional case, with
$\alpha=1, \alpha^{\vee}=2$, the
graph of a function $g:[0,1]\to \mathbb R$ as well as those of
$I,-I $ and $f=\pa g$ where
    $I(s)=\inf _{0\leq u\leq s}g(u)$.
   \hfill $\diamondsuit$.
$$
\begin{array}{cc}\includegraphics[scale=0.3]{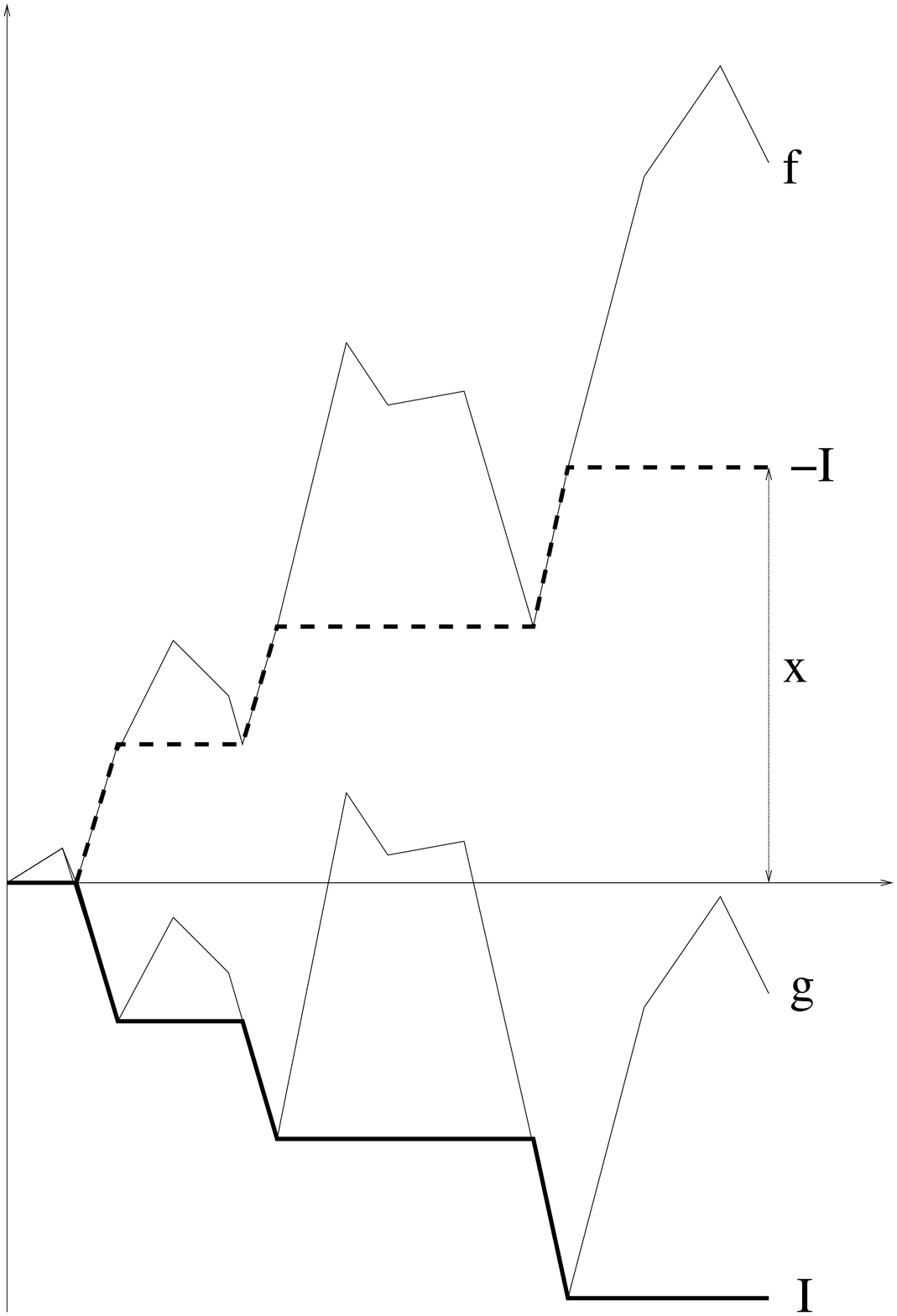}\\
\text{Fig. 1}
\end{array}
$$
\subsection{Relation with Littelmann path operators.}\label{littelpit}
Using Proposition \ref{pit} ({\it iv}) we can define generalized
Littelmann
transformations. Recall that Littelmann operators are defined
    on  paths
with values in the dual space $\mathfrak a^*$ of
    some real Lie algebra $\mathfrak a$. The image of
a path is either another path or the symbol $\bf 0$ (actually the zero
element in
the $\mathbb Z$-module generated by all paths). We define continuous
versions of
these operators.
\begin{definition}
Let $\pi:[0,T]\to V$ be a continuous path satisfying $\pi(0)=0$, and  
$x\in \mathbb R$, then
$E_{\alpha}^x\pi$ is the unique path such that $$\pa E_{\alpha}^x\pi=\pa  
\pi\quad
\text{and}\quad\alpha^{\vee}(E_{\alpha}^x\pi(T))=\alpha^{\vee}(\pi(T))+x 
$$
if $
    -2\alpha^{\vee}(\pi(T))+2\inf_{0\leq t\leq T}\alpha^{\vee}(
\pi(t))\leq x\leq -2\inf_{0\leq t\leq T}\alpha^{\vee}(
\pi(t))$ and $E_{\alpha}^x\pi=\bf 0$ otherwise.
\end{definition}
One checks easily that
    $E_{\alpha}^0\pi=\pi$ and
    $E_{\alpha}^x E_{\alpha}^y\pi=E_{\alpha}^{x+y}\pi$ as long as
    $E_{\alpha}^y\pi\ne \bf 0$.
When $\alpha$ is a root and $\alpha^{\vee}$ its coroot, in some root
system,
    then  $E_{\alpha}^2$ and $E_{\alpha}^{-2}$
coincide with the Littelmann operators $e_{\alpha}$ and $f_{\alpha}$,
    defined in
\cite{littel}. 
Recall that a path $\pi$ is called integral if its endpoint $\pi(T)$
 is in the weight
lattice and, for each simple root $\alpha$,
 the
minimum of the function $\alpha^{\vee}(\pi(t))$ is an integer. The class of 
integral paths
is invariant under the Littelmann operators.
For such paths, the action of a Pitman
transform can be expressed through Littelmann operators by
\begin{equation}\label{pit/lit}
\pa \pi=e_{\alpha}^{n_\alpha}(\pi)
\end{equation}
where $n_\alpha$ is the largest integer
$n$ such that $e_{\alpha}^{n}(\pi)\ne \bf 0$.

    \subsection{Braid relations}
An important property of the Pitman transforms is the following result.
\begin{theorem}\label{braid}
Let $\alpha,\beta\in V$ and   $\alpha^{\vee},\beta^{\vee}\in V^{\vee}$
     be  such that
$\alpha^{\vee}(\alpha)=\beta^{\vee}(\beta)=2$, and
$\alpha^{\vee}(\beta)< 0,
\beta^{\vee}(\alpha)< 0$ and
$\alpha^{\vee}(\beta)
\beta^{\vee}(\alpha)=4\cos^2\frac{\pi}{n}$, where $n\geq 2$
    is some  integer,
then one has
    $$\pa\pb\pa \ldots =\pb\pa\pb\ldots$$ where there are $n$
factors in each product.
\end{theorem}
We shall prove Theorem \ref{braid} as a corollary to the result of
section
\ref{papb}. Note that if $\alpha^{\vee}(\beta)=
\beta^{\vee}(\alpha)=0$ then $\pa\pb=\pb\pa$ by a simple computation.  
For crystallographic angles (i.e. $n=2,3,4,6$)
a proof of Theorem \ref{braid}
 could also be deduced
from Littelmann's theory (see \cite{littel2} or  \cite{kash}).
We shall provide still another (hopefully
more conceptual)
proof for these angles in section \ref{rep-pw}, see Remark \ref{remprf}.
    The
general case seems to be new.

\subsection{A formula for $\pa\pb\pa\pb\ldots$}\label{papb}
Let $\alpha,\beta\in V$ and   $\alpha^{\vee},\beta^{\vee}\in V^{\vee}$
be such
that $\alpha^{\vee}(\beta)< 0$ and $\beta^{\vee}(\alpha)< 0$. By
Proposition \ref{pit} {(\it i)} we can - and will - assume by rescaling that
$\alpha^{\vee}(\beta)=\beta^{\vee}(\alpha)$, without changing $\pa$ and
$\pb$.
We use the notations
$$
\rho=-\frac{1}{2}\alpha^{\vee}(\beta)=
-\frac{1}{2}\beta^{\vee}(\alpha),\quad
X(s)=\alpha^{\vee}( \pi(s)),\quad Y(s)=\beta^{\vee} (\pi(s)).
$$
\begin{theorem}\label{formula} Let $n$ be a positive integer, if
$\rho\geq\cos\frac{\pi}{n}$, then one has
\begin{eqnarray}\label{forpapb}
(\underbrace{\pa\pb\pa\ldots}_{ \text{$n$ terms}} )
\pi(t)&=&\pi(t)-\inf_{t\geq s_0\geq s_1\geq \ldots\geq s_{n-1}\geq
0}\bigl(\sum_{i=0}^{n-1}T_i(\rho)Z^{(i)}(s_i)\bigr)\alpha\nonumber\\&&
-\inf_{t\geq s_0\geq s_1\geq \ldots\geq s_{n-2}\geq
0}\bigl(\sum_{i=0}^{n-2}T_i(\rho)Z^{(i+1)}(s_i)\bigr)\beta
\end{eqnarray}
where $Z^{(k)}=X$ if $k$ is even and $Z^{(k)}=Y$ if $k$ is odd. The
$T_k(x)$ are the Tchebycheff polynomials defined by
$T_0(x)=1$, $T_1(x)=2x$, and
$2xT_k(x)=T_{k-1}(x)+T_{k+1}(x)$ for $k\geq 1$.
\end{theorem}
    The Tchebycheff polynomials satisfy
$T_k(\cos\theta)=\frac{\sin(k+1)\theta}{\sin\theta}$ and, in particular,
under the assumptions on $\rho$ and $n$, one has
$T_k(\rho)\geq 0$ for all $k\leq n-1$.

Assuming Theorem \ref{formula} we obtain Theorem \ref{braid}.

\medskip
{ Proof of Theorem \ref{braid}.}
Let $\alpha^{\vee}(\beta)=
\beta^{\vee}(\alpha)=-2\cos\frac{\pi}{n}$, then one has
$T_{n-1}(\rho)=0$
    and the last term in the coefficient of $\alpha$ in the right hand
side of
    (\ref{forpapb})
    vanishes.
    It follows by inspection
     that this term equals the coefficient of $\alpha$ in
    the analogous formula for $\underbrace{\pb\pa\pb\ldots}_{ \text{$n$
    terms}}\pi(t)$. A similar argument works for the coefficient of
$\beta$.
   \hfill $\diamondsuit$

    \medskip
    The proof of Theorem \ref{formula} will be by
    induction on $n$. It is easy to check the formula for $n=1$ or $2$.
    We shall do the induction in
    sections 2.5 and 2.6.
    

\subsection{Two intermediate lemmas}\label{IL}
\begin{lemma}\label{aux}
Let $X:[0,t]\to \mathbb R$ be a continuous functions with $X(0)=0$ and
let
$t_0=\sup\{s\geq 0\,|\, X_s=\inf_{s\geq u\geq 0}X_u\}$, then for all
$u\leq  t_0$
one has  $$\inf_{t\geq s\geq u}(X(s)-2\inf_{s\geq w\geq
0}X(w))=-\inf_{u\geq v\geq 0}X(v).$$
\end{lemma}
{  Proof.} This is obtained as a byproduct of the proof in section
\ref{prfpit}. Again it is perhaps more convincing to stare at Fig. 1
than to
give a formal proof.
\hfill$\diamondsuit$

\medskip
Elaborating on this we obtain the next result.
\begin{lemma}\label{fund}
Let $X$ and $Y$ be continuous functions, such that $X(0)=Y(0)=0$, then
\begin{eqnarray*}
&\inf_{t\geq s\geq 0}\bigl (X(s)+\inf_{s\geq u\geq 0} Y(u)\bigr)
=
\inf _{t\geq s\geq 0}X(s)+\\&
\inf_{t\geq s\geq 0}\Bigl(X(s)-2\inf _{s\geq u\geq 0}X(u)
+\inf_{s\geq u\geq 0} \bigl(Y(u)+\inf_{u\geq v\geq 0}
X(v)\bigr)\Bigr)\nonumber
\end{eqnarray*}
\end{lemma}
{ Proof.}
The first term is $I=\inf_{t\geq s\geq u\geq 0}(X(s)+Y(u))$.
Let $t_0$ be, as in Lemma \ref{aux},
    the last time when  $X$ reaches its minimum over $[0,t]$, then
$$I=\inf\bigl(\inf_{t_0\geq u\geq 0} (Y(u)+X(t_0));
\inf_{t\geq s\geq u\geq t_0\geq 0}
(Y(u)+X(s))\bigr)$$
Let $J$ be the second term in the identity to be proved, then $$J=
\inf_{t\geq s\geq 0}\bigl[X(s)-2\inf _{s\geq u\geq 0}X(u)
+\inf_{s\geq u\geq 0} (Y(u)+\inf_{u\geq v\geq 0} X(v))\bigr]
+X(t_0)$$
    Introduce again the time $t_0$, then
\begin{eqnarray*}
J&=&
\inf_{t\geq s\geq u\geq 0}\bigl(X(s)-2\inf _{s\geq w\geq 0}X(w)
+Y(u)+\inf_{u\geq v\geq 0} X(v)\bigr)
+X(t_0)\\&=&
\inf\Bigl(\inf_{\substack{t\geq s\geq u\geq 0\\ t_0\geq u}}\bigl(Y(u)+
X(s)-2\inf _{s\geq w\geq 0}X(w)
+\inf_{u\geq v\geq 0} X(v)+X(t_0)\bigr);\\&&\qquad
\inf_{t\geq s\geq u\geq t_0}
\bigl(Y(u)+
X(s)-2\inf _{s\geq w\geq 0}X(w)
+\inf_{u\geq v\geq 0} X(v)+X(t_0)\bigr)\Bigr)
\end{eqnarray*}
but if $u\leq t_0$ then by lemma \ref{aux}
    one has $\inf_{t\geq s\geq u}(X(s)-2\inf_{s\geq w\geq
0}X(w))=-\inf_{u\geq v\geq 0}X(v)$. If  $t_0\leq u$ then
$\inf _{s\geq w\geq 0}X(w)=X(t_0)$,
therefore
\begin{eqnarray*}
J&=&
\inf\Bigl(\inf_{t_0\geq u\geq 0}\bigl(
Y(u)+X(t_0)\bigr);
\inf_{0\geq u\geq t_0}
\bigl(Y(u)+\inf_{t\geq s\geq u}
X(s)\bigr)\Bigr)\\&=&
\inf\bigl(\inf_{t_0\geq u\geq 0} (Y(u)+X(t_0));
\inf_{t\geq s\geq u\geq t_0\geq 0}
(Y(u)+X(s))\bigr)\\
&=&I.
\end{eqnarray*}

\subsection {End of proof of Theorem \ref{formula}}
Assume the result of the Theorem  holds for some  $n$ with $n$ even.
Then
$\underbrace{\pa\pb\pa\ldots}_{\text{$n+1$ terms}}=
\underbrace{\pa\pb\pa\ldots}_{ \text{$n$ terms}}\pa$, and one has
\begin{eqnarray*}
\alpha^{\vee}(\pa \pi(s))&=&X(s)-2\inf_{s\geq u\geq 0}X(u)\\
\beta^{\vee}(\pa \pi(s))&=&Y(s)+2\rho\inf_{s\geq u\geq 0}X(u)
\end{eqnarray*}
therefore, by induction hypothesis
\begin{eqnarray*}
\underbrace{\pa\pb\pa\ldots}_{ \text{$n+1$ terms}}
\pi(t)&=&\underbrace{\pa\pb\pa\ldots}_{ \text{$n$ terms}}
(\pa \pi)(t)
\\&=&\pi(t)-\inf_{t\geq s\geq 0}X(s)\alpha
-\inf_{t\geq s_0\geq s_1\geq \ldots\geq  s_{n-1}\geq
0}\Bigl(\sum_{i=0}^{n-1}\hat Z^{(i)}(s_i)\Bigr)\alpha
\\
&&
-\inf_{t\geq s_0\geq s_1\geq \ldots\geq s_{n-2}\geq
0}\Bigl( \sum_{i=0}^{n-2}\hat Z^{(i+1)}(s_i)\Bigr)\beta
\end{eqnarray*}
where $$\hat Z^{(i)}_{\alpha}(s)=\left\{
\begin{matrix}X(s)-2\inf_{s\geq u\geq 0}X(u)&\text{ for $i$ even} \\
Y(s)+2\rho\inf_{s\geq u\geq 0}X(u)&\text{ for $i$
odd.}\end{matrix}\right.$$
The coefficient of $\alpha$ in the above expression has the form
\begin{eqnarray*}
&&H_{\alpha}=-\inf_{t\geq s\geq 0}T_0(\rho)X(s)\\&&
-\inf_{t\geq s\geq
0}\Bigl(T_0(\rho)X(s)-2\inf_{s\geq u\geq 0}T_0(\rho)X(u)+
\inf_{s\geq u\geq 0}\bigl(\Gamma(u)+\inf_{u\geq v\geq 0}T_0(\rho)
X(v)\bigr)\Bigr)
\end{eqnarray*}
where
\begin{eqnarray*}
\Gamma(u)&=&T_1(\rho) Y(u)+2\rho T_1(\rho)\inf_{u\geq v\geq 0}X(v)+
\\&&\quad\inf_{u\geq
u_2\geq u_3\geq\ldots\geq
u_{n-1}\geq 0}\Bigl(\sum_{i=2}^{n-1}T_i(\rho)\hat Z^{(i)}(u_i)\Bigr)
-T_0(\rho)\inf_{u\geq v\geq 0}X(v)\\
&=&T_1(\rho) Y(u)+T_2(\rho)\inf_{u\geq v\geq 0}X(v)+
\\&&\quad\inf_{u\geq
u_2\geq u_3\geq\ldots\geq
u_{n-1}\geq 0}\Bigl(\sum_{i=2}^{n-1}T_i(\rho)\hat Z^{(i)}\Bigr)
\end{eqnarray*}
so that we can apply lemma \ref{fund} to transform it into
$$H_{\alpha}=-\inf_{t\geq s\geq 0}\bigl(T_0(\rho)X(s)
+\inf_{s\geq u\geq 0}\Gamma(u)\bigr)$$

Let us prove by induction on $k$ that
\begin{eqnarray*}
H_{\alpha}&=&-\inf_{t\geq u_0\geq u_1\geq\ldots\geq u_{2k}}\Bigl(
\sum_{i=0}^{2k}T_i(\rho)Z^{(i)}(u_i)+
    W_{k}(u_{2k-1})\biggr)
\end{eqnarray*}
with
\begin{eqnarray*}
W_{k}(v)&=&
\inf_{v\geq u_{2k}\geq u_{2k+1}\geq \ldots\geq u_{n-1}\geq
0}\Bigl(\sum_{i=2k}^{n-1}T_{i}(\rho)(\hat Z^{(i)}(u_i)\Bigr)
\end{eqnarray*}
Indeed the formula holds for $k=1$ by the computation above.
Assume this holds for some $k$ then one  has
\begin{eqnarray*}
&&H_{\alpha}=-\inf_{t\geq u_0\geq u_1\geq\ldots\geq u_{2k}}\Bigl(
\sum_{i=0}^{2k}T_i(\rho)Z^{(i)}(u_i)+
    W_{k}(u_{2k-1})\biggr)\\&&\qquad=
-\inf_{t\geq u_1\geq u_2\geq\ldots\geq u_{2k-1}}\biggl(
\sum_{i=0}^{2k-1}T_i(\rho)Z^{(i)}(u_i)+\inf_{u_{2k-1}\geq v\geq
0}T_{2k}(\rho)X(v)+\\
&&\qquad\inf_{u_{2k-1}\geq v\geq 0}
    \Bigl((T_{2k}(\rho) X(v)
-2\inf_{v\geq w\geq 0}T_{2k}(\rho)X(w))+\\&&\qquad\inf_{w\geq z\geq
0}\bigl(R_k(z)+\inf_{z\geq \tau\geq
0}T_{2k}(\rho)X(\tau)\bigr)\Bigr)\biggr)
\end{eqnarray*}
where
\begin{eqnarray*}
R_k(z)&=&T_{2k+1}(\rho)Y(z)+2\rho T_{2k+1}(\rho)\inf_{z\geq \tau\geq
0}X(\tau)
+\\&&\quad\inf_{z\geq u_{2k+2}\geq \ldots u_{n-1}}
\bigl(\sum_{i=2k+2}^{n-1}\hat Z^{(i)}(u_{i})\bigr)-
\inf_{z\geq \tau\geq 0}T_{2k}(\rho)X(\tau)\\&=&
T_{2k+1}(\rho)Y(z)+T_{2k+2}(\rho)\inf_{z\geq \tau\geq 0}X(\tau)
+\inf_{z\geq u_{2k+2}\geq \ldots u_{n-1}}
\bigl(\sum_{i=2k+2}^{n-1}\hat Z^{(i)}(u_{i})\bigr)
\end{eqnarray*}
where we used $2\rho T_{2k+1}(\rho)-T_{2k}(\rho)=T_{2k+2}(\rho)$.
    Applying Lemma 2.1. we get
    \begin{eqnarray*}
H_{\alpha}&=&-\inf_{t\geq u_1\geq u_2\geq\ldots\geq u_{2k-1}}\Bigl(
\sum_{i=0}^{2k-1}T_i(\rho)Z^{(i)}(u_i)+\\&&\qquad\inf_{u_{2k-1}\geq
v\geq 0}\bigl(T_{2k}(\rho)X(v)+
\inf_{w\geq z\geq
0}R_k(z)\bigr)\Bigr)\\&=&
-\inf_{t\geq u_0\geq u_1\geq\ldots\geq u_{2k+2}}\biggl(
\sum_{i=0}^{2k+2}T_i(\rho)Z^{(i)}(u_i)+W_{k+1}(u_{2k+1})\biggr)
\end{eqnarray*}
Taking $k=n$ gives the required formula for $H_{\alpha}$.
For the coefficient of $\beta$, remark that
$$
\underbrace{\pa\pb\pa\ldots}_{ \text{$n+1$ terms}}
\pi(t)=\pa(\underbrace{\pb\pa\pb\ldots}_{ \text{$n$ terms}}
\pi)(t)
$$
and the formula for $n+1$ follows immediately from the formula at step
$n$ for
$\underbrace{\pb\pa\pb\ldots}_{ \text{$n$ terms}}$.
    The case where $n$ is odd is
treated in a similar way.
\hfill$\diamondsuit$
\subsection{Pitman transformations for Coxeter and Weyl  
groups}\label{sec_cox}
Let $W$ be a Coxeter group, i.e. $W$ is generated by a
 finite set $S$ of reflections of a  
real vector space $V$, and $(W,S)$ is a
Coxeter system (see
\cite{bo}, \cite{humphreys}). For each $s \in S$, let $\alpha_s\in V$  
and  $\alpha_s^\vee \in V^\vee$, where $V^\vee$ is the dual space of
$V$,  such that $s=s_{\alpha_s}$  is the reflection  
associated to $\alpha_s$ (see (\ref{salpha})).
Then $\alpha_s$ is called the simple root associated with $s\in S$ and  
$\alpha_s^\vee$ its coroot.

Denote by
$\PP_s$ the Pitman transform associated with  the pair  
$(\alpha_s,\alpha^\vee_s)$.
By the results of the preceding sections, 
the $\PP_s;s\in S$ form a representation of the monoid generated by idempotents
satisfying the braid relations. Such a monoid occurs in the theory of Hecke
algebras for $q=0$, and in the calculus of Borel orbits (see e.g. \cite{knop}
where this monoid is called Richardson-Springer monoid).

Let     $H_s$ be the closed
half space $H_s=\{v\in V|\alpha^\vee_s(v)\geq 0\}$.
Let $w\in W$ and let $w=s_1\ldots s_l$ be a reduced decomposition of
$w$, where $l=l(w)$ is the length of $w$.
By Theorem \ref{braid} and
a fundamental result of Matsumoto (\cite {bo} Ch. IV, $n^o$ 1.5,
    Proposition 5)  the operator
$\PP_{s_1}\ldots \PP_{s_l}$ depends only on $w$, and not on the chosen
reduced
decomposition. We shall denote by $\PP_w$ this operator.
\begin{prop}
Let $w\in W$, $L_w=\{s\in S\,|\,l(sw)<l(w)\}, R_w=\{s\in  
S\,|\,l(ws)<l(w)\}$.
For any path $\pi$, the path  $\PP_w\pi$ lies in the convex cone
$\cap_{s\in L_w}H_s$,  one has $\PP_s\PP_w=\PP_w$ for all $s\in L_w$  
and $\PP_w\PP_s=\PP_w$ for all $s\in R_w$.
\end{prop}
{ Proof.} If $l(sw)<l(w)$ then $w$ has a reduced
decomposition $w=ss_{1}\ldots s_{k}$ therefore
$\PP_w=\PP_s\PP_{s_{1}}\ldots
\PP_{s_{k}}$ and $\PP_w\pi=\PP_s(\PP_{s_{1}}\ldots
\PP_{s_{k}}\pi)$ lies in $H_s$ by Proposition \ref{pit} ({\it ii}).  
Furthermore one has $\PP_s\PP_w=\PP_w$ since $\PP_s$ is an involution (see
Proposition \ref{pit} ({\it ii}) ). Similarly $\PP_w\PP_s=\PP_w$ when
$l(ws)<l(w)$.
 \hfill$\diamondsuit$

\begin{cor} If $W$ is finite and $w_0$ is the longest element, then
$\PP_{w_0}\pi$ takes values in the closed Weyl chamber $\overline   
C=\cap_{s\in S}H_s$, furthermore
    $\PP_{w_0}$ is an idempotent and  
$\PP_w\PP_{w_0}=\PP_{w_0}\PP_w=\PP_{w_0}$ for all $w \in W$.
\end{cor}
    Assume now that $W$ is a finite Weyl group, associated with a
     weight lattice
    in $V$.
    Recall that paths taking values in the Weyl chamber $\overline C$  
are called
    dominant paths in \cite{littel}, and that
    the set $B\pi$ of all (nonzero) paths obtained by
    applying products of  Littelmann operators to a dominant path $\pi$
    is called the Littelmann module. 
     From the connection between Pitman's and
    Littelmann's operators,
    given in section \ref{littelpit}, we deduce the following (see also  
\cite{littel2}).
\begin{cor}\label{pit-littel}
Let $\pi$ be a dominant integral path, then 
 a path $\eta$ belongs to the Littelmann module $B\pi$ if and only
if $\eta$ is integral and 
$\pi=\PP_{w_0}\eta$. 
\end{cor}
Indeed for any path $\eta$ and $x$ such that $E_\alpha^x\eta\ne\bf 0$
one has $\PP_{\alpha}E_\alpha^x\eta=\PP_{\alpha}\eta$, therefore
$\PP_{w_0}E_\alpha^x\eta=\PP_{w_0}\PP_{\alpha}E_\alpha^x\eta=\PP_{w_0}\eta$.
It follows that the set of paths whose image by $\PP_{w_0}$ is $\pi$ is stable under
the action of Littelmann operators. If $\eta$ is  an integral path
such that $\PP_{w_0}\eta=\pi$, and $w_0=s_1\ldots s_n$ is a reduced
decomposition, then by section \ref{littelpit} the sequence
$\eta, \PP_{\alpha_n}\eta,\PP_{\alpha_{n-1}}\PP_{\alpha_n}\eta,\ldots,\pi$ is
obtained by successive applications of Littelmann operators
 therefore they all belong to the
Littelmann module $B\pi$.  \hfill$\diamondsuit$

\medskip

Let us come back to the general case of a finite Coxeter group.
We shall now study the set of all paths $\eta$
such that $\PP_{w}\eta$ is a given dominant
path. Let $w=s_{1}\ldots s_{q}$ be  a reduced decomposition.
Let $\eta$ be a path such that $\eta(0)=0$ and
     $\pi=\PP_{w}\eta$ is a dominant path. Denote
$\eta_0=\pi,\eta_q=\eta$, and
    $\eta_j=\PP_{s_{j+1}}\ldots \PP_{s_q}\eta_q$ for $j=1,2,\ldots ,q-1$,
then by
Proposition \ref{pit} ({\it iv}) for all $j=1,2,\ldots ,q$ the
    path $\eta_{j}$ is uniquely specified among paths $\gamma$ such
    that $\PP_{s_j}\gamma=\eta_{j-1}$, by the number $x_j=-
    \inf_{0\leq t\leq T}\alpha_{s_j}^{\vee}( \eta_j(t))\in
    [0,\alpha_{s_j}^{\vee}( \eta_{j-1}(T))]$. It follows that
$\eta=\eta_q$ is uniquely
    specified, among all
    paths $\gamma$ such that $\PP_{w_0}\gamma=\pi$ by the sequence
    $x_1,x_2,\ldots, x_q$.
    These coordinates are   subject to the inequalities
    $0\leq x_j\leq \alpha_{s_j}^{\vee}(\eta_{j-1}(T))$.
    From
    $$
\eta_{j-1}(T)=\eta_{j}(T)
    +x_{j}\alpha_{s_j}$$
    one obtains  
$$\pi(T)=\eta_0(T)=\eta_j(T)+\sum_{l=1}^jx_l\alpha_{s_l}$$
    therefore the inequality $0\leq x_j\leq
\alpha_{s_j}^{\vee}(\eta_{j-1}(T))$
     reads
    $$0\leq x_j\leq
     \alpha_{s_j}^{\vee}(\pi(T))-\sum_{l=1}^{j-
1}x_l\alpha_{s_j}^{\vee}(\alpha_{s_l}).$$
It follows that
    the set of all
paths $\eta$ such that $\PP_{w}\eta=\pi$ can be parametrized
by a subset of  the convex polytope
$$K_\pi=\{(x_1,\ldots, x_q)\in \mathbb R^q|\,0\leq x_j\leq
\alpha_{s_j}^{\vee}(
    \pi(T))-\sum_{l=1}^{j-1}x_{l}
\alpha_{s_j}^{\vee}(\alpha_{s_l});j=1,\ldots, q\}.$$
    The path $\eta$ corresponding to the point $(x_1,\ldots, x_q)$ is
specified by
    the equalities
    $$
\eta_{j-1}(T)=\eta_{j}(T)
    +x_{j}\alpha_{s_j}$$ where
    $\eta_j=\PP_{s_{j+1}}\ldots \PP_{s_q}\eta$.
In the case of a Weyl group, it follows from  \cite{littel2}
that the subset of $K_\pi$ corresponding to paths $\eta$ such that
$\mathcal P_w
\eta=\pi$ is the intersection of $K_{\pi}$ with a certain  convex cone
which
does not depend on $\pi$. This convex cone is quite difficult to
describe, see \cite{beze2}. Also
we do not know if a similar result holds for all finite Coxeter
groups.
We hope to come back to these questions in future work.
\section{A representation theoretic formula for $\PP_{w}$}\label{rep-pw}
\subsection{Semisimple groups}\label{sec_semisimple} We recall some  
standard terminology. We
consider
a simply connected complex semisimple Lie group $G$, associated with a  
root system $R$.
Let $H$ be a maximal torus, and $B^+,B^-$ be corresponding
opposite Borel subgroups with unipotent
radicals $N^+,N^-$.
Let $\alpha_i,i\in I,$ and $\alpha_i^{\vee},i\in I,$ be the
simple positive roots and coroots,  and $s_i$ the
corresponding reflections
in the Weyl group $W$.
Let $e_i,f_i, h_i,i\in I,$ be  Chevalley
generators of the Lie algebra of $G$. 
One can  choose
representatives $\overline w\in G$ for  $w\in W$ by putting
$\overline{s_i}=\exp(-e_i)\exp(f_i)\exp(-e_i)$ and
$\overline{vw}=\overline v\,\overline w$ if
$l(v)+l(w)=l(vw)$ (see \cite{foze} (1.8), (1.9)).  The Lie algebra of  
$H$, denoted by
     $\mathfrak h$ has a Cartan decomposition $\mathfrak
    h=\mathfrak a+i\mathfrak a$ such that the roots
     $\alpha_i$ take real values on the real vector space $\mathfrak
    a$. Thus $\mathfrak a$ is generated by $\alpha_i^{\vee}, i\in I$ and  
its dual $\mathfrak a^*$ by $\alpha_i, i\in I$.
The set of weights  is the lattice $P=\{\lambda \in \mathfrak a^*;  
\lambda(\alpha_i^\vee)\in \mathbb Z, i \in I\}$ and the set of dominant
weights is $P^+=\{\lambda \in \mathfrak a^*; \lambda(\alpha_i^\vee)\in  
\mathbb N, i \in I\}$.
  For each  $\lambda\in P^+$, choose a representation space
$V_{\lambda}$
with a highest weight vector $v_{\lambda}$,
and an invariant inner product on $V_{\lambda}$ for which $v_{\lambda}$
is
a unit vector.

\begin{lemma} \label{positive}
For any dominant weight $\lambda$, $w\in W$
 and indices $i_1,\ldots,i_n\in I$
    one has
$$\langle e_{i_1}\ldots e_{i_n}\overline
wv_{\lambda},v_{\lambda}\rangle\geq 0$$
\end{lemma}
Proof. This is an immediate consequence  of Lemma 7.4 in \cite{beze2}.
$\diamondsuit$

\medskip

Let
$(\omega_i,i\in I)\in P^I$ be the fundamental weights, characterized by the  
relations $\omega_i(\alpha_j^\vee)=\delta_{i,j}, j \in
I$. The principal minor associated with $\omega_i$ is the function on $G$ given
by
$$\Delta^{\omega_i}(g)=\langle gv_{\omega_i},v_{\omega_i}\rangle$$
see \cite{beze} and \cite{foze}.
If $g\in G$ has a Gauss decomposition $g=[g]_-[g]_0[g]_+$ with
$[g]_-\in N^-,\, [g]_0\in H,\, [g]_+\in N^+$, then one has
\begin{equation}\label{gauss}
\Delta^{\omega_i}(g)=[g]_0^{\omega_i}=e^{\omega_i(\log[g]_0)}.
\end{equation}

\subsection {Some auxiliary path transformations}
    We shall now introduce some  path transformations.
    \begin{definition} Let $n_i:[0,T]\to \mathbb R^+, i \in I$,
     be a family of strictly positive continuous functions,
     and let $a:(0,T]\to  \mathfrak a$ be
     a  continuous map such that $$\int_{0^+}e^{-\alpha_i(
    a(s))}n_{i}(s)ds<\infty $$ we define, for $0 < t \leq T$,
    $$\T_{i,n}a(t)=a(t)+\log\left (\int_0^te^{-\alpha_i(
    a(s))}n_{i}(s)ds\right)\alpha^{\vee}_i.$$
    \end{definition}
    Observe that in general the maps $t\mapsto a(t)$ and  
$t\mapsto\T_{i,n}a(t)$ need not be
continuous at 0. For all that follows,  consideration of the case 
  $n_i\equiv
1$ in the above definition would be sufficient for our purposes, but the proofs
would be the same as the general
case.

    Let $R^{\vee}$ be the root system dual to $R$, namely the roots of
    $R^{\vee}$ are the coroots of $R$ and vice versa,
    and denote by
    $\PP_{\alpha_i^{\vee}},i\in I,$ the corresponding Pitman
transformations on $\mathfrak a$. Let $\pi$ be a
continuous path in $\mathfrak a$, with $\pi(0)=0$. For $\varepsilon >0$,
let $D_{\varepsilon}$ be the dilation operator  
$D_{\varepsilon}\pi(t)=\varepsilon \pi(t)$. A simple
application of Laplace method yields the following
\begin{equation}\label{laplace}
\mathcal P_{\alpha_i^{\vee}}\pi=\lim_{\varepsilon\to
0}D_{\varepsilon}\T_{i,n}D_{\varepsilon}^{-1}\pi .
\end{equation}
We shall establish, in section \ref{repT},
    a representation theoretic formula for a product
$\T_{i_k,n}\ldots \T_{i_1,n}$ corresponding to a minimal decomposition
$w=s_{i_1}\ldots s_{i_k}$ in the Weyl group. Using this formula we shall
use (\ref{laplace}) to get a formula for the Pitman transform.
\subsection{ A group theoretic interpretation of the operators
$\T_{i,n}$}

Let $a$ be a smooth path in $\mathfrak a$
    and let $b$ be the path in the Borel subgroup $B^+=HN^+$ solution to
the differential equation
$$\frac{d}{dt}b(t)=\left(\frac{d}{dt}a(t)+\sum_{i\in I}n_i(t)e_i\right)b(t);
\qquad b(0)=id.$$
The following expression is easy to check.
\begin{lemma}\label{sol}
\begin{eqnarray}\label{bt}
&&\qquad\qquad\qquad\qquad b(t)=e^{a(t)}+
e^{a(t)}\sum_{k\geq 1}\sum_{i_1,\ldots ,i_k\in I^k}
\\&&\left(\int_{t\geq t_1\geq t_2\geq\ldots
\geq t_k\geq 0}e^{-\alpha_{i_1}(a(t_1))}n_{i_1}(t_1)\ldots
e^{-\alpha_{i_k}(a(t_k))}n_{i_k}(t_k)dt_1\ldots
dt_k\right)e_{i_1}\ldots e_{i_k}\nonumber
\end{eqnarray}
\end{lemma}
Observe that this  expression is well 
 defined in each finite dimensional representation of  
$G$ since  the operators $e_i$ are
nilpotent and  this sum has only a finite number of nonzero terms.
It is always in this context that we shall use this formula.
\begin{lemma}\label{dec_b}
For any $t>0$ and $w\in W$ one has
$$\Delta^{\omega_i}(b(t)\overline w)>0.$$
\end{lemma}
{Proof.} By  eq.\  (\ref{bt})
    one has
    \begin{eqnarray}\label{bt2}
   &&  \Delta^{\omega_i}(b(t)\overline w)=\langle e^{a(t)}\overline w
v_{\omega_i},v_{\omega_i}\rangle+\\&&\,
\qquad\sum_{r\geq 1}\sum_{i_1,\ldots ,i_r\in I^r}
\int_{t\geq t_1\geq t_2\geq\ldots
\geq t_r\geq 0}\langle e^{a(t)}
e^{-\alpha_{i_1}(a(t_1))}n_{i_1}(t_1)\ldots\nonumber \\ &&\qquad\quad
\ldots
e^{-\alpha_{i_r}(a(t_r))}n_{i_r}(t_r) e_{i_1}\ldots e_{i_r}\overline
wv_{\omega_i},v_{\omega_i}\rangle\,dt_1\ldots
dt_r \nonumber
\end{eqnarray}
    which
    is a sum of nonegative terms by
Lemma \ref{positive}. Furthermore, since $v_{\omega_i}$ is a
    highest weight vector, there exists some sequence
    $i_1,\ldots ,i_r$, such that $e_{i_1}\ldots e_{i_r}\overline
wv_{\omega_i}$
    is a nonzero multiple of $v_{\omega_i}$, and the $n_i$ do not vanish,
     therefore the sum is positive. \hfill$\diamondsuit$

\medskip
It follows in particular
that, according to the terminology of \cite{foze},
    $b(t)$  belongs to the double Bruhat cell $ B_+\cap B_-w_0B_-$,
     and that $b(t)\overline w$ has a Gauss
decomposition $b(t)\overline w=[b(t)\overline w]_-[b(t)\overline
w]_0[b(t)\overline w]_+$
    for all $t>0$.

Now comes the main result of this section.
    \begin{theorem}\label{main}
    Let  $w\in W$ and
$w=s_{i_1}\ldots s_{i_k}$ be a reduced decomposition, then
    the $H$ part in the Gauss decomposition
    of $b(t)\overline w$ is equal to
$$\exp(\T_{i_k,n}\ldots \T_{i_1,n}a(t)).$$
\end{theorem}
The fact that the path $\T_{i_k,n}\ldots \T_{i_1,n}a(t)$ is well
defined is part
of the Theorem.
By the uniqueness of the
Gauss decomposition the preceding result implies
\begin{cor}
The path
$$\T_{i_k,n}\ldots \T_{i_1,n}a(t)$$ depends only on $w$ and $n$ and not  
on the
chosen reduced
decomposition of $w$.
\end{cor}We shall denote by $\T_{w}a$ the resulting path (it depends on  
$n$). We thus have
\begin{equation}\label{eq_bt}
[b(t)\overline w]_0=e ^{\T_wa(t)}.
\end{equation}

\medskip
{ Proof of Theorem \ref{main}.}
    The proof  is  by induction on the length of
$w$.  Let $s_i$ be such that $l(ws_i)=l(w)+1$. We assume that the $H$
part of
the Gauss decomposition of
$b(t)\overline w$ is $\T_{i_k,n}\ldots \T_{i_1,n}a(t)$ as required.
By (\ref{gauss}) it is then enough to prove that for all $t>0$ and
$i,j\in I$
    one has
$$\Delta^{\omega_j}(b(t)\overline w\overline
s_i)=\Delta^{\omega_i}(b(t)\overline w)$$
if $i\ne j$ and
$$\Delta^{\omega_i}(b(t)\overline w\overline
s_i)=\Delta^{\omega_i}(b(t)\overline w)
\int_0^te^{-\alpha_i(
    \T_wa(s))}n_{i}(s)ds.$$
    The claim for $i\ne j$ follows from Proposition 2.3 in \cite{foze},  
it
remains
    to check the case $i=j$.
    \begin{lemma}\label{zero}
$$\frac{\Delta^{\omega_i}(b(t)\overline w\overline s_i)}
{\Delta^{\omega_i}(b(t)\overline w)}\to
_{t\to 0}0.$$
\end{lemma}
{Proof. }
 From the decomposition (\ref{bt2}), the fact that all
    terms are positive and that the $n_i$ are positive continuous
functions,
we see that as $t\to 0$ one has $\Delta^{\omega_i}(b(t)\overline w)\sim
c_1t^{l_1}$
and $\Delta^{\omega_i}(b(t)\overline w\overline s_i)\sim
c_2t^{l_2}$ for some $c_1,c_2>0$, where $l_1$ (resp.\  $l_2$) is the
number of terms in the decomposition of $\omega_i-w(\omega_i)$ (resp.\
$\omega_i-ws_i(\omega_i)$) as a sum of simple roots. Since
$l(ws_i)>l(w)$ the weight $w(\omega_i)-ws_i(\omega_i)$ is positive, and  
one
has $l_2>l_1$. \hfill$\diamondsuit$
\begin{lemma}
Let $w=s_{i_1}\ldots s_{i_k}$ be a reduced decomposition, and let
$b^{w}(t)=[b(t)\overline w]_0[b(t)\overline w]_+$, then
     one has
$$\frac{d}{dt}b^{w}(t)=\left(\frac{d}{dt}\T_{i_k,n}\ldots \T_{i_1,n}
a(t)+\sum_{j\in I}n_j(t)e_j\right)b^w(t).$$
\end{lemma}
{Proof.}
We do this by induction on the length of $w$.
Assume this is true for $w$ and let $s_i$ be such that
$l(ws_i)=l(w)+1$, then
one has
$$\frac{d}{dt}b^w(t)=\left(\frac{d}{dt}\T_wa(t)+\sum_jn_je_j\right)b^w(t
)$$
    therefore
$$\frac{d}{dt}b^w(t)\overline
s_i=\left(\frac{d}{dt}\T_wa(t)+\sum_jn_j(t)e_j\right)b^w(t)
\overline s_i$$
Since $b^w(t)\in B^+$, by \cite{beze}, \cite{foze},
the Gauss decomposition of $b^w(t)\overline s_i$ has the form
$$b^w(t)\overline s_i=\exp(\beta(t)f_i)b^{ws_i}(t)$$ with $\beta(t)>0$
for $t>0$,
    and one has, since $f_i$
commutes with all   $e_j$ for $j\ne i$.
\begin{eqnarray*}
\frac{d}{dt}b^{ws_i}(t)&=&\frac{d}{dt}\left[\exp(-
\beta(t)f_i)b^{w}(t)\overline
s_i\right]\\&=&-\left(\frac{d}{dt}\beta(t)\right)f_i\exp(-
\beta(t)f_i)b^{w}(t)\overline
s_i+\\&&\qquad\qquad
\exp(-
\beta(t)f_i)\left(\frac{d}{dt}\T_wa(t)+\sum_{j}n_j(t)e_j\right)b^{w}(t)\
\overline
s_i\\&=&-\frac{d}{dt}\beta(t)f_ib^{ws_i}(t)+\\&&\quad
\left(\frac{d}{dt}\T_wa(t)+\sum_{j}n_j(t)e_j+n_i(t)\beta(t)h_i+n_i(t)
\beta^2(t)f_i\right)b^{ws_i}(t)
\\&=&
\Biggl[\left(\frac{d}{dt}\beta(t)+\frac{d}{dt}\alpha_i(\T_wa(t))+n_i(t)\beta^2(t)  \right)f_i
+\\&&\qquad
\frac{d}{dt}\T_wa(t)+n_i(t)\beta(t)h_i+\sum_{j}n_j(t)e_j\Biggr]b^{ws_i}(
t)
\end{eqnarray*}
Since $b^{ws_i}(t)\in B_+$, one has
    $\frac{d}{dt}\beta(t)+\frac{d}{dt}
    \alpha_i(T_wa(t))+\beta^2(t) =0$ therefore
$$\beta(t)=\frac{e^{-\alpha_i(\T_wa(t))}}{C+\int_0^te^{-
\alpha_i(T_wa(s))}
    n_i(s)ds}$$
    for some constant $C\geq 0$. Integrating the $H$ part of the Gauss
    decomposition of $b^{ws_i}(t)$
     we see that this part is equal to
     \begin{equation}\label{hpart}
     \exp(\T_wa(t))
     \exp(C'+\log(C+\int_0^te^{-\alpha_i(T_wa(s))}n_i(s)ds))h_i
     \end{equation}
      therefore
     $$\frac{\Delta^{\omega_i}(b(t)\overline w\overline
s_i)}{\Delta^{\omega_i}(b(t)\overline w)}=
     \exp(C')(C+\int_0^te^{-\alpha_i(\T_wa(s))}n_i(s)ds)$$
     and $C=0$ by Lemma \ref{zero}.
      We conclude that
$$\beta(t)=\frac{e^{-\alpha_i(\T_wa(t))}}{\int_0^te^{-\alpha_i(T_wa(s))}
    n_i(s)ds}.$$ This implies that
     \begin{eqnarray*}
\frac{d}{dt}b^{ws_i}(t)
&=&
\left[\frac{d}{dt}\T_wa(t)+n_i(t)\frac{e^{-
\alpha_i(\T_wa(t))}}{\int_0^te^{-\alpha_i(T_wa(s))}
    n_i(s)ds}h_i+\sum_{j}n_j(t)e_j\right]b^{ws_i}(t)
\\&=&\left[\frac{d}{dt}\T_{i,n}\T_{w}a(t)+\sum_{j}n_j(t)e_j\right]b^{ws_
i}(t)
\end{eqnarray*}
as required.
\hfill$\diamondsuit$

\medskip

    From (\ref{hpart}) we obtain
     $$\frac{\Delta^{\omega_i}(b(t)\overline w\overline
s_i)}{\Delta^{\omega_i}(b(t)\overline w)}
     =\Delta^{\omega_i}(e^{-\T_wa(t)}b^w(t)\overline s_i)
     =\exp(C')\int_0^te^{-\alpha_i(T_wa(s))}n_i(s)ds$$
     Differentiating with respect to $t$ we get
     \begin{eqnarray*}
     \frac{d}{dt}e^{-\T_wa(t)}b^w(t)\overline
     s_i&=&e^{-\T_wa(t)}\sum_jn_j(t)e_je^{\T_wa(t)}
     e^{-\T_wa(t)}b^w(t)\overline
     s_i\\
     &=&\left(\sum_ie^{-\alpha_j(\T_wa(t))}n_j(t)e_j\right)
     e^{-\T_wa(t)}b^w(t)\overline
     s_i\\
     \end{eqnarray*}
     where $e^{-\T_wa(t)}b^w(t)\in N$.
     It follows that
     \begin{eqnarray*}\frac{d}{dt}
     \left\{\frac{\Delta^{\omega_i}(b(t)\overline w\overline s_i)}
     {\Delta^{\omega_i}(b(t)\overline
     w)}\right\}
     &=&\left\langle\left(\sum_j
     e^{-\alpha_j(\T_wa(t))}n_j(t)e_j\right)
     e^{-\T_wa(t)}b^w(t)\overline
     s_iv_{\omega_i},v_{\omega_i}\right\rangle\\
     &=&e^{-\alpha_i(\T_wa(t))}n_i(t)\langle
     e_i\overline s_iv_{\omega_i},v_{\omega_i}\rangle\\
     &=&e^{-\alpha_i(\T_wa(t))}n_i(t)
     \end{eqnarray*}
     therefore $C'=0$. This proves the claim for $i=j$ and finishes the
proof of
     Theorem \ref{main}.
  \hfill   $\diamondsuit$
    \begin{cor}
    The transformations $\T_{i,n}$ satisfy the braid relations,
    $$\underbrace{\T_{i,n}\T_{j,n}\ldots}_{ \text{$m(i,j)$ terms}}
     =\underbrace{\T_{j,n}\T_{i,n}\ldots}_{ \text{$m(i,j)$ terms}}$$
where $m(i,j)$ is the Cartan integer $\alpha_i(\alpha_j^\vee)$.
    \end{cor}
    \begin{remark}\label{remprf}
    In the case of rank two groups, the braid relations of the above
corollary
    and an application of Laplace method yield the braid relations for
Pitman
    operators as in Theorem \ref{braid}, in the case of cristallographic
angles
    ($\pi/m,m=2,3,4,6)$. It is instructive to give an elementary  
derivation
of the
    braid relations for the $\T_{i,n}$ in the simplest nontrivial case
namely type
    $A_2$ (i.e., $m=3$). In this case the relations amount to
   \begin{equation}\label{naf-gtr}
\int_0^tds\int_0^sdr\ F(r)\frac{G(s)}{G(r)}\frac{H(t)}{H(s)}
= \int_0^tds\int_0^sdr\ F(r)\frac{\tilde G(s)}{\tilde G(r)}\frac{\tilde
H(t)}{\tilde H(s)} ,
\end{equation}
for some positive continuous functions $F,G,H$,  where
$$\tilde G(s) = \left( \int_0^s  G(r) H(r)^{-1} dr \right)^{-1} G(s)$$
and
$$\tilde H(s) = \left( \int_0^s  G(r) H(r)^{-1} dr \right) H(s).$$
   This can be checked directly by an application of Fubini's theorem,
or an
    integration by parts. Similar but more complicated formulas  
correspond to
    the other crystallographic angles $\pi/4$ and $\pi/6$.

 From (3.7) one recovers, by the method of Laplace, the identity
\begin{equation}\label{naf-tr}
x\tr (z\ut y)\tr (y\tr z) = (x\tr y)\tr z ,
\end{equation}
for continuous functions $x,y,z$ with $x(0)=y(0)=z(0)=0$
and (non-associative) binary operations $\ut$ and $\tr$ defined by
\begin{equation}
(x\tr y)(t)=\inf_{0\le s\le t} [x(s)-y(s)+y(t)] ,
\end{equation}
\begin{equation}
(x\ut y)(t)=\sup_{0\le s\le t} [x(s)-y(s)+y(t)] .
\end{equation}
This is equivalent to the $n=3$ braid relation for the Pitman 
transforms.
For a `queueing-theoretic' proof, which some readers might find 
illuminating,
see~\cite{oc}.  Lemma~\ref{fund} is a special case.

    \end{remark}
    \subsection{Representation theoretic formula for  
$\PP_{w}$}\label{repT}
    Let $w\in W$, and let $\lambda$ be a dominant weight,
    then $\lambda- w\lambda$ can be decomposed as a linear
    combination of
    simple positive roots $\lambda-w\lambda=\sum_{i\in I}u_i\alpha_i$
    where $u_i$ are nonnegative integers. If $(j_1,\ldots,j_r)\in I^r$ is
a    sequence such that
     $\langle e_{j_1}\ldots
e_{j_r}\overline w v_{\lambda},v_{\lambda}\rangle\ne
0$, then the number of $k$'s in the sequence $j_1,\ldots,j_r$
is equal to $u_k$. In particular
the number $r$ depends only on $w$ and $\lambda$.
    We let
    $S(\lambda,w)$ denote the set of sequences $(j_1,\ldots,j_r)\in I^r$
    such that
$\langle e_{j_1}\ldots
e_{j_r}\overline w v_{\lambda},v_{\lambda}\rangle\ne
0$.
    Using (\ref{bt2}) and  (\ref{eq_bt})
    we  obtain the following  expression
    \begin{prop}\label{T-rep}
    Let $a$ be a path in $\mathfrak a$, and $\lambda$ a dominant weight,
    then one has
    \begin{eqnarray*}
    &&\langle
e^{\T_{w}a(t)}v_{\lambda},v_{\lambda}\rangle
    = e^{\lambda(a(t))}\sum_{(j_1,\ldots,j_r)\in S(\lambda,w)}
    \int_{t\geq t_1\geq \ldots \geq t_r\geq 0}\\&&\qquad e^{-\alpha_{j_1}
    (a(t_1))-\ldots -\alpha_{j_r}( a(t_r))}
    n_{j_1}(t_1)\ldots n_{j_r}(t_r)dt_1\ldots dt_r\langle e_{j_1}\ldots
    e_{j_r}\overline w v_{\lambda},v_{ \lambda}\rangle
    \end{eqnarray*}
    \end{prop}
    Let $w\in W$ and let
    $\PP^{\vee}_{w}$ denote the Pitman transformation on $\mathfrak a$ for the  
dual root
system
    $R^{\vee}$,
    by (\ref{laplace}), one has $$\mathcal P^{\vee}_{w}\pi=\lim_{\varepsilon\to
0}D_{\varepsilon}\T_{w}D_{\varepsilon}^{-1}\pi.$$ Using Laplace method,   
Lemma \ref{positive} and
     Proposition \ref{T-rep} applied to fundamental weights, we now  
obtain
    the following expression for the Pitman transform (notice that $W$  
acts on $\mathfrak a^*$ and on $\mathfrak a$ by duality).
    \begin{theorem} \label{bjoc}
    (Representation theoretic formula for the Pitman transforms).
    
   Let $w\in W$, for each path $\pi$ on $\mathfrak a$, one has
    \begin{equation}\label{formidable}
    \PP^{\vee}_{w}\pi(t)=\pi(t)-\sum_{i\in I}
    \inf_{\substack{j_1,\ldots, j_r\in S(\omega_i,w)\\
    t\geq t_1\geq t_2\ldots \geq t_r\geq
       0}}\left(\alpha_{j_1}(\pi(t_1))+\ldots+\alpha_{j_r}(\pi(t_r))\right)\alpha_i^{\vee}
       \end{equation}
    \end{theorem}
    This formula can be seen as a generalization of the formula in
    Theorem \ref{formula}.
    Observe that  sequences $j_1,\ldots j_r$ such as the ones occuring in
the
    theorem have appeared already in \cite{beze2} under
    the name of {\bf i}-trails. It is interesting to note that such sequences
    appear here naturally by an application of the Laplace method (sometimes
    called "tropicalization" in the algebraic litterature).
    
     By Corollary 1, we see that Theorem \ref{bjoc} provides a
representation
     theoretic formula for the dominant path in some Littelmann module, which is
     independent of any choice of a reduced decomposition of $w_0$.
\begin{remark}\label{conj}
As noted before, formula \ref{formidable} has a similar structure as
  formula
\ref{forpapb} (when $\rho=\cos\frac{\pi}{n}$). We conjecture that such formulas
exist for arbitrary Coxeter groups, i.e. for $w\in W$ there exists $r$ and 
a set $S(s,w)\subset S^r$  such that
\begin{equation}
    \PP^{}_{w}\pi(t)=\pi(t)-\sum_{s\in S}
    \inf_{\substack{s_1,\ldots, s_r\in S(s,w)\\
    t\geq t_1\geq t_2\ldots \geq t_r\geq
      0}}\left(\alpha^{\vee}_{s_1}(\pi(t_1))+\ldots+\alpha^{\vee}_{s_r}(\pi(t_r))\right)\alpha_s.
       \end{equation}
       However we do not know how to interpret these sets $S(s,w)$.
      \end{remark}

\section{Duality}
  \subsection{An involution on dominant paths} As in section  
\ref{sec_cox}, we consider  a Coxeter system $(W,S)$  generated
by a  set $S$ of reflections of
$V$. We assume now that the group $W$  is finite and let $w_0$ be the longest
element. We fix some
$T >0$ and for any continuous path $\pi:[0,T] \to V$ such that  
$\pi(0)=0$ we let
 $$\kappa \pi(t)=\pi(T-t)-\pi(T).$$ Clearly for all paths $\kappa^2\pi=\pi$.
  We will show that the transformation $I=\PP_{w_0}\kappa(-w_0)$
is an involution on the set of dominant paths, which generalizes the  
Sch\"utzenberger involution (see section
\ref{schutz} for the connection).
\subsection{Codominant paths and co-Pitman operators}
A path $\pi$ is called $\alpha$-dominant if 
$\alpha^{\vee}(\pi(t))\geq 0$ for all $t$.
It is called $\alpha$-codominant if $\kappa\pi$ is $\alpha$-dominant
or, in other words, if
$\alpha^{\vee}(\pi(t))\geq \alpha^{\vee}(\pi(T))$ for all $t$.
Finally it is called codominant if it is $\alpha$-codominant for all $\alpha$.
Let us define the co-Pitman operators 
$\ea=\kappa\pa \kappa$, given by the formula
$$\ea\pi(t)=\pi(t)-\inf_{t\leq s\leq T}\alpha^{\vee}(\pi(s))\alpha+
\inf_{0\leq s\leq T}\alpha^{\vee}(\pi(s))\alpha$$
One checks the following

$$\pa\kappa\pa=\pa,\quad \ea^2=\ea,\quad \ea\pa=\ea,\quad\pa\ea=\pa$$
Furthermore for all paths $\pi$ one has
$$\ea\pi(T)=s_{\alpha}\pa\pi(T).$$
A few properties of $\ea$ are gathered in the following lemma, whose proof is
left to the reader.

\begin{lemma}
({\it i}) $\ea\pi$ is the unique path $\eta$ satisfying 
$\eta(T)=s_{\alpha}\pa\pi(T)$ and $\pa\eta=\pa\pi$.

({\it ii})
$\ea\pi$ is the unique path $\eta$ such that $\pa\eta=\pa\pi$ and 
$\eta$ is $\alpha$-codominant.

({\it iii}) If $\pi$ is $\alpha$-dominant, then $\ea\pi$ is the unique path such that
$\pa\eta=\pi$ and $\eta(T)=s_{\alpha}(\pi(T))$.

({\it iv}) $\ea\pi=\pi$ if and only if $\pi$ is $\alpha$-codominant.
\end{lemma}
The transformations $\ea$ play the same role with respect to the Littelmann
operators $f_{\alpha}$ as the transformation $\pa$ with respect to $e_{\alpha}$
(see (\ref{pit/lit})). 
\begin{lemma}
The $\ea$ satisfy the braid relations.
\end{lemma}
Proof.   Follows from $\ea=\kappa\pa\kappa$, $\kappa^2=id$ and 
 the braid relations for the $\pa$.\hfill $\diamondsuit$
 
\medskip

One can therefore define $\ew$ for $w\in W$, and
$\ewo=\ewo^2$ is a projection onto the set of codominant paths. Furthermore 
for all $w\in W$ one has
$$\ew=\kappa{\mathcal P}_{w}\kappa$$
 In particular
 $$\ewo=\kappa{\mathcal P}_{w_0}\kappa.$$
\subsection{An endpoint property}
In this section we prove the following  result, which is crucial for
applications to Brownian motion.

\begin{prop}\label{e-duality}For any path $\pi$ one has
$$\ewo\pi(T)=w_0\pwo\pi(T)$$ 
\end{prop}
Since 
 $\pwo\ewo=\pwo$ 
it is enough to check this identity for $\pi$ a codominant  path
(or for a dominant path using $\ewo\pwo=\ewo$).
\begin{lemma}Let $\pi$ be a codominant path, 
let  $w\in W$ and  $\alpha$ be such that
$l(s_{\alpha}w)>l(w)$, then $\pw\pi$ is $\alpha$-codominant.

\end{lemma}

Proof. First we check the result for dihedral groups. With the  
notations of
\ref{papb}, let $\pi$ be a
$\alpha$- and $\beta$-codominant path, and let $n$ be such that
$\rho>\cos\frac{\pi}{n}$,
  then one has $\alpha^{\vee}(\pi(T))\leq\alpha^{\vee}(\pi(t))$ and
  $\beta^{\vee}(\pi(T))\leq\beta^{\vee}(\pi(t))$ for all $t\leq T$.
   It follows that in the
  computation of $\underbrace{\pb\pa\pb\ldots}_{ \text{$n$ terms}}\pi(T)$
  using formula (\ref{forpapb})
  the infimum is obtained for $s_0=s_1=\ldots=T$, therefore
  (assuming $n$ odd for definiteness)
  \begin{eqnarray*}
  &&\alpha^{\vee}(\pb\pa\pb\ldots  \pi(T))
  =\\&&
  \quad\alpha^{\vee}(\pi(T))+2\rho\bigl[\beta^{\vee}(\pi(T))
+T_1(\rho)
  \alpha^{\vee}(\pi(T))+\ldots
  +T_{n-1}(\rho)\beta^{\vee}(\pi(T)))\bigr]\\&&
  \quad-2\bigl[\alpha^{\vee}(\pi(T))+T_1(\rho)
  \beta^{\vee}(\pi(T))+\ldots
  +T_{n-2}(\rho)\beta^{\vee}(\pi(T)))\bigr]\\&&
  = T_{n-1}(\rho)\alpha^{\vee}(\pi(T))+T_{n}(\rho)\beta^{\vee}(\pi(T))
  \end{eqnarray*}
  where we have used the recursion relation of the $T_k$.
   On the other hand, for $t\leq T$ one has
  \begin{eqnarray*}
&&\alpha^{\vee}(\pb\pa\pb\ldots
\pi(t))=\alpha^{\vee}(\pi(t))+\\
&&\quad
2\rho\inf_{t\geq s_0\geq \ldots\geq s_{n-1}\geq
0}\bigl[\beta^{\vee}(\pi(s_0))+T_1(\rho)\alpha^{\vee}(\pi(s_1))+\ldots
T_{n-1}(\rho)\beta^{\vee}(\pi(s_{n-1})\bigr]\\
&&\quad
-2\inf_{t\geq s_0\geq \ldots\geq s_{n-2}\geq
0}\bigl[\alpha^{\vee}(\pi(s_0))+T_1(\rho)\beta^{\vee}(\pi(s_1))+\ldots
T_{n-2}(\rho)\beta^{\vee}(\pi(s_{n-2})\bigr]
\end{eqnarray*}
In this expression let us replace, inside the
$\inf_{t\geq s_0\geq s_1\geq \ldots\geq s_{n-1}\geq
0}$ each $2\rho T_k(\rho)$ by $T_{k-1}(\rho)+T_{k+1}(\rho)$. We obtain
\begin{eqnarray*}
&&\inf_{t\geq s_0\geq s_1\geq \ldots\geq s_{n-1}\geq 0}
        \bigl[2\rho\beta^{\vee}(\pi(s_0))+
        (T_0(\rho)+T_2(\rho))\alpha^{\vee}(\pi(s_1))+\ldots\\ &&\qquad
        (T_{n-2}(\rho)+T_n(\rho))\beta^{\vee}(\pi(s_{n-1}))\bigr]\geq
        \\
        &&\quad
        \inf_{\substack{t\geq s_0\geq s_1\geq \ldots\geq s_{n-1}\geq 0\\
        t\geq u_1\geq\ldots\geq u_{n-1}\geq 0}}
        \bigl[2\rho\beta^{\vee}(\pi(s_0))+T_2(\rho)\alpha^{\vee}(\pi(s_1))
        +\ldots\\&&\qquad
         + T_n(\rho)\beta^{\vee}(\pi(s_{n-1}))+T_0(\rho)\alpha^{\vee}(\pi(u_1))
        +\ldots+ T_{n-2}(\rho)\beta^{\vee}(\pi(u_{n-1})\bigr]=
        \\
        &&\quad
        \inf_{t\geq s_0\geq \ldots\geq s_{n-1}\geq 0}
        \bigl[2\rho\beta^{\vee}(\pi(s_0))+T_2(\rho)\alpha^{\vee}(\pi(s_1))
        +\ldots +T_n(\rho)\beta^{\vee}(\pi(s_{n-1}))\bigr]
        +
        \\
        &&
        \qquad\inf_{t\geq s_0\geq \ldots\geq s_{n-2}\geq 0}
         \bigl[\alpha^{\vee}(\pi(s_0))+T_1(\rho)\beta^{\vee}(\pi(s_1))+\ldots
        T_{n-2}(\rho)\beta^{\vee}(\pi(s_{n-2}))\bigr]
\end{eqnarray*}
Furthermore
\begin{eqnarray*}
&&\alpha^{\vee}(\pi(t))+
\\ && \quad\inf_{t\geq s_0\geq \ldots\geq s_{n-1}\geq 0}
        \bigl[2\rho\beta^{\vee}(\pi(s_0))+
        T_2(\rho)\alpha^{\vee}(\pi(s_1))+\ldots+
        T_{n}(\rho)\beta^{\vee}(\pi(s_{n-1})\bigr]\geq
        \\
        &&
        \quad\inf_{t\geq s_0\geq \ldots\geq s_{n-2}\geq 0}
         \bigl[\alpha^{\vee}(\pi(s_0))+T_1(\rho)\beta^{\vee}(\pi(s_1))+\ldots
+ T_{n-2}(\rho)\beta^{\vee}(\pi(s_{n-2})\bigr]
\\ && \qquad+T_{n-1}(\rho)\alpha^{\vee}(\pi(T))+
T_n(\rho)\beta^{\vee}(\pi(T))
\end{eqnarray*}
Putting everything together we obtain
$$\alpha^{\vee}(\pb\pa\pb\ldots
\pi(t))\geq \alpha^{\vee}(\pb\pa\pb\ldots
\pi(T))$$
and $\pb\pa\pb\ldots \pi$ is $\alpha$-codominant.
The case of $n$ even is similar. This proves the claim for dihedral  
groups.

Consider now a general Coxeter system.
We do the  proof  by induction on $l(w)$.
The claim is true if $l(w)=0$. If
it is true for some $w$, let
$s_\beta\in S$ be such that $l(s_{\beta}w)>l(w)$. Let now
$\alpha$ be such that
$l(s_{\alpha}s_{\beta}w)>l(s_{\beta}w)>l(w)$.
Let $n$ be the order of $s_{\alpha}s_{\beta}$, and
 $$w=s_{\alpha}w_1=s_{\alpha}s_{\beta}w_2
=s_{\alpha}s_{\beta}s_{\alpha}w_3=\ldots=s_{\alpha}s_{\beta}\ldots
w_k.$$
where  $k$ is the smallest integer such that $$l(w)>l(w_1)>\ldots  
 >l(w_k)\quad\hbox{
and }\quad l(s_{\alpha}w_k)>l(w_k),\,
l(s_{\beta}w_k)>l(w_k).$$ Since $l(s_\alpha s_\beta w)=l(w_k)+k+2$ one has
$k+2\leq n$.
By induction hypothesis,
${\mathcal P}_{w_k}(\pi)$ is both $\alpha$ and $\beta$ codominant.
Then it follows from the dihedral case that
$\PP_{s_{\beta}}\pw=\pb\pw\pi=\pb\pa\pb\ldots \mathcal{P}_{w_k}\pi$ is $\alpha$-codominant.
$\diamondsuit$

\begin{lemma}Let $\pi$ be a codominant path and
 $w\in W$, then  $\pw\pi$ is the unique path $\eta$ such that
${\mathcal E}_{w^{-1}}\eta=\pi$, and $w(\pi(T))=\eta(T)$.
\end{lemma}
The proof is by induction on $l(w)$, using the preceding lemma. Let
$l(s_{\alpha}w)=l(w)+1$, then $\pw \pi$ is $\alpha$-codominant, therefore
$\pa\pw\pi$ is the unique path $\eta$
such that $\ea\eta=\pw\pi$, and $\eta(T)=s_{\alpha}\pw\pi(T)$.
\hfill $\diamondsuit$

\medskip

Proposition \ref{e-duality} is the special case $w=w_0$ in the last lemma.

\begin{lemma} \label{Lem_dual} $(-w_0)\PP_{w_0}=\PP_{w_0}(-w_0)$.
\end{lemma}
{ Proof.}
If $\alpha$ is a simple root, then $\tilde\alpha=-w_0\alpha$ is also a  
simple root and
$\tilde\alpha^\vee=-\alpha^\vee w_0$. It follows easily that
$(-w_0)\PP_\alpha (-w_0)=\PP_{\tilde\alpha}.$

If $w_0=\alpha_1\cdots\alpha_r$ is a reduced expression we thus have
$$\PP_{w_0}(-w_0)=\PP_{\alpha_1}\cdots
{\PP_{\alpha_r}}(-w_0)=(-w_0)\PP_{\tilde\alpha_1}\cdots
{\PP_{\tilde\alpha_r}}=(-w_0)\PP_{w_0}$$
since $w_0=\tilde\alpha_1\cdots \tilde\alpha_r$. \hfill$\diamondsuit$

\begin{theorem}\label{Theo_dual} The transformation $I=\PP_{w_0}(-w_0)\kappa$
 has the following
properties:

({\it i}) $I^2=\PP_{w_0}$;

({\it ii}) The  restriction of $I$ to dominant paths is an involution;

({\it iii}) $I\PP_{w_0}=I$;

({\it iv}) (Duality relation) For all paths $\pi$,  
one has $$I\pi(T)=\PP_{w_0}\pi(T);$$ in particular, one has
$I\pi(T)=\pi(T)$ when $\pi$ is dominant.
\end{theorem}
{ Proof.} By  Lemma  
\ref{Lem_dual},
$$I^2=\PP_{w_0}\kappa(-w_0)\PP_{w_0}(-w_0)\kappa=\PP_{w_0}\kappa\PP_{w_0}\kappa=\PP_{w_0}
\ewo=\PP_{w_0}$$
this proves (i) and implies (ii) since $\PP_{w_0}\pi=\pi$ when $\pi$ is  
dominant. This also give
$$I\PP_{w_0}=I^3=I^2I=I$$
since the image by $I$ of any path is dominant. Finally 
$I=\PP_{w_0}\kappa(-w_0)=\kappa\ewo(-w_0)$, and  Proposition \ref{e-duality}
gives ({\it iv}). \hfill $\diamondsuit$

\medskip

Property ({\it iv}) will be important for the first
 proof of the Brownian motion
property.
\subsection{Symmetry of a Littlewood-Richardson construction}\label{LR}
The concatenation $\pi\star\eta$ of two paths $\pi: [0,T] \to V$ $\eta:  
[0,T] \to V$ is defined in Littelmann \cite{littel} as the
path
$\pi\star\eta:[0,T]\to V$ given by
$\pi\star\eta(t)=\pi(2t)$, when
$0
\leq t
\leq T/2$ and $\pi\star\eta(t)=\pi(T)+\eta(2(t-T/2))$ when $T/2 \leq t  
\leq T$.

\begin{lemma} \label{lem_concat}For all $w\in W$ one has
$\PP_w(\pi\star\eta)=\PP_w(\pi)\star\eta'$, where  
$\PP_{w_0}(\eta')=\PP_{w_0}(\eta)$.
\end{lemma}
{Proof.} One uses induction on
the length $l(w)$ of $w$. When $l(w)=1$  it is easy to see that  
$\PP_w(\pi\star\eta)=\PP_w(\pi)\star\eta'$ where
$\PP_{w}(\eta)=\PP_w(\eta')$. Since $\PP_{w_0}\PP_w=\PP_{w_0}$ the  
claim is thus true in this case. Suppose that it holds for
elements of length
$n$. Let
$w=w_1s$ where
$l(w)=n+1,l(w_1)=n$, then one has
$$\PP_{w}(\pi\star\eta)=\PP_{w_1}\PP_{s}(\pi\star\eta)=\PP_{w_1}(\PP_{ 
s}(\pi)\star \eta')$$
  where $\PP_{w_0}\eta'=\PP_{w_0}\eta$. Now by induction hypothesis
$$\PP_{w_1}(\PP_{s}(\pi)\star \eta')=(\PP_{w_1}\PP_{s})(\pi)\star  
\eta''$$
where $\PP_{w_0}\eta''=\PP_{w_0}\eta'$, and therefore   
$\PP_{w_0}\eta''=\PP_{w_0}\eta$.\hfill$\diamondsuit$
\medskip

In the case of Weyl groups, Littelmann has given the
following analogue of the Littlewood-Richardson construction: Let
$\pi$ and
$\eta$ be two integral dominant paths defined on $[0,T]$, then the set
$$LR(\pi,\eta)=\{\pi\star\mu\,|\, \mu\in B\eta, \pi\star\mu \mbox{  
is dominant}\} $$
gives a parametrization of the decomposition into irreducible representations
of the  tensor product of the representations  
with highest weights $\pi(T)$ and $\eta(T)$.
 By Theorem \ref{Theo_dual} ({\it iii}), one has 
$I(\eta)(T)=\eta(T)$ and $I(\pi)(T)=\pi(T)$, therefore
$LR(I(\eta),I(\pi))$ gives a parametrization of the decomposition 
of the tensor product of  
the representations with highest weights
$\eta(T)$ and $\pi(T)$.

\begin{prop} 
The map
$I:LR(\pi,\eta))\to LR(I(\eta),I(\pi))$
is a bijective involution, which preserves the end points.
\end{prop}
{ Proof.} Let  $\pi\star\mu \in LR(\pi,\eta)$. By Lemma  
\ref{lem_concat} there is a path $\xi$ such that
$$I(\pi\star\mu)=\PP_{w_0}(\kappa(-w_0)(\pi
\star\mu))=\PP_{w_0}(\kappa(-w_0)(\mu) \star \kappa(-w_0)(\pi))
=\PP_{w_0}(\kappa(-w_0)(\mu))\star  
\xi $$  and
$\PP_{w_0}\xi=\PP_{w_0}(\kappa(-w_0) (\pi))=I(\pi)$.  By {\it(iii)} of Theorem  
\ref{Theo_dual}, one has
$I(\mu)=I(\eta)$ thus $I(\pi\star\eta)\in  
LR(I(\eta),I(\pi))$. One checks easily that $I$ preserves integrality,
and the other properties follow from 
Theorem \ref{Theo_dual}.\hfill$\diamondsuit$

\subsection{Connection with the  Sch\"utzenberger involution} \label{schutz}
In the case of a Weyl group of type $A_{d-1}$
  the transform $\PP_{w_0}$ is connected with the  
Robinson, Schensted and Knuth (RSK) correspondence : Let us consider  
a word
$v_1v_2\cdots v_n$ written with the alphabet
$\{1,2,\cdots,d\}$. Let $(P(n),Q(n))$ be the pair of  tableaux  
associated  with this word by  RSK with column insertion (see, e.g.,
\cite{ful}). Let
$ \mathfrak a=\{(x_1,\cdots,x_d)\in  \mathbb R^d; \sum_{i=1}^d x_i=0\}$  
and let
$(e_i)$ be the image in $\mathfrak a$ of the canonical basis of $  
\mathbb R^d$.  We identify $v_i$
with the path $\eta_i:t\mapsto te_{v_i}, 0 \leq t \leq 1,$ and we  
consider the path
$\pi=\eta_1
\star
\eta_2
\cdots \star
\eta_n$. Then $\PP_{w_0}\pi$ is the path obtained by taking  the  
successive shapes of $Q(1),Q(2),\cdots,Q(n)$ (see Littelmann
\cite{littel},\cite{littel3}, or \cite{oc} for a connection with queuing theory).
  Let us
consider the pair
$(\tilde P(n),\tilde Q(n))$ associated by the RSK algorithm to the word
$v_n^*\cdots v_1^*$ where
$v^*=d+1-v$.
The Sch\"utzenberger involution is the map which associates the tableau  
$\tilde Q(n)$ to the tableau $Q(n)$
(see \cite{fomin}, \cite{ful}, \cite{van Leeuw}). The path  associated  
with the word $v_n^*\cdots v_1^*$ is $I(\pi)$. Thus  $I$
is a generalization of this involution. Note that $I$ makes sense not only for
Weyl groups, but for any finite Coxeter group.
\section{Representation of Brownian motion in a Weyl chamber}
\subsection{Brownian motion in a Weyl chamber}\label{BMWC}
In this section we recall some basic facts about Brownian motion in Weyl
chambers.

We consider  a Coxeter system $(W,S)$ generated by a set $S$ of  
reflections
of an  euclidean space $V$ and we assume that $W$ is finite.  We  
shall denote
by $C$ the interior of a fundamental domain for the action of $W$ on $V$
(a  Weyl chamber), and by $\overline C$ its closure.

If $W$ is the Weyl group of a complex semi-simple Lie algebra  
$\mathfrak g$, with
compact form
$\mathfrak g_{\mathbb R}$, then $V$ is identified
with $\mathfrak a^*$, the dual space of the Lie algebra  of a maximal  
torus $T$,  and
the Weyl chamber $\overline C=\mathfrak {\overline a}^*_+$
   can be identified with the orbit space of
$\mathfrak g_{\mathbb R}^*$
under the coadjoint action of the simply connected compact group $K$  
with Lie algebra $\mathfrak g_{\mathbb R}$
(up to some  identification of the
   walls). Let $Z$ be a Brownian motion with values in $\mathfrak
g_{\mathbb R}^*$,
whose covariance is the Killing form.  It is well known that the
image  of  $Z$ in the quotient space $ \mathfrak g_{\mathbb R}^*/K$  
remains in the interior
of the
Weyl chamber for all times $t>0$, even if the starting point is inside
some
wall.   Since the transition
probabilities of $Z$
are invariant under the coadjoint action it follows that this
image,  under the quotient map, is a Markov
   process on $\overline C$.  A description of this
Markov process can be done in terms of Doob's conditionning, namely the
process
is obtained from a  Brownian motion $X$ on $V=\mathfrak a^*$, killed  
at the boundary
of the Weyl chamber, by means of a Doob transform with
   respect to the function
$$h(v)=\prod_{\alpha\in R^+}\alpha^\vee(v), \ v \in V,$$ (where $R^+$  
is the set of
positive roots)
which is the unique, up to a scaling
factor, positive harmonic function on $\overline C$ which vanishes on
the
boundary (see \cite{bi3}). Recall that, by the reflection principle, the
transition probabilities for the Brownian motion killed at the boundary
of the Weyl chamber are
\begin{equation}\label{p0}
p^{0}_t(x,y)dy=\sum_{w\in W}\varepsilon(w)p_t(x,wy)dy, \ x,y \in  
\overline C,
\end{equation}
where $p_t(x,y)dy$
   are the transition probabilities for Brownian motion $X$,
given by the  Gaussian kernel on $\mathfrak a^*$ whose covariance is
that of the Brownian motion.
Thus the probability transitions for the Doob's process are
\begin{equation}\label{q-p}
q_t(x,y)dy=\frac{h(y)}{h(x)}\sum_{w\in W}\varepsilon(w)p_t(x,wy)dy
\end{equation}
for $x\in C$.
   These probability transitions can be continued by continuity to
    $x\in \overline C$ , in particular to $x=0$.

   For a general finite Coxeter group, formula
   (\ref{p0})  still gives the probability transitions of Brownian motion killed
   at the boundary of the Weyl chamber.
    Let $h$ be the product of the  
positive coroots, defined as the
   linear forms corresponding to the hyperplanes
   of the reflections in the group $W$, taking the signs so that they are
positive
   inside the Weyl chamber, then
   the function $h$ is still the only (up to a multiplicative constant) positive
   harmonic function vanishing on the boundary, and the equation (\ref{q-p})
   defines the semi-group of what we call the Brownian motion in the  
fundamental chamber
   $\overline C$ of $V$. 
   
   We shall prove that the
    Pitman operator $\PP_{w_0}$ applied to Brownian motion in
   $V$ yields a Brownian motion in the Weyl chamber. We shall give two very
   different proofs of
   this. The first one uses in an essential way the duality relation of
   Proposition 
   \ref{e-duality}  and a classical result in queuing theory.
  The second one uses a
   random walk approximation and relies on
   Littelmann theory and Weyl's character formula. It is  valid only for
   Weyl groups. We have chosen to present this second proof because it 
   emphazises the close connection between Brownian paths and Littelmann paths.
   \subsection{Brownian motion with a drift}
   We now consider a Brownian motion in $V$ with invariant covariance,  
but with
a drift
   $\xi\in C$. Its transition probabilities  are
now
   $$p_{t,\xi}(x,y)=p_t(x,y)
   \exp(\langle \xi,y-x\rangle-\frac{\Vert \xi\Vert^2}{2}t)$$
   Actually the distribution of this Brownian motion on the  
$\sigma$-field
   $\mathcal F_t$ generated by the coordinate functions $X_s, s \leq t,$  
  on
the
   canonical space,
   is absolutely continuous with respect to the one of the centered
   Brownian motion, with density
   $$\exp(\langle \xi,X_t-X_0\rangle-\frac{\Vert \xi\Vert ^2}{2}t). $$
   Consider such a Brownian motion in $V$ with drift $\xi$, starting  
inside the
chamber
   at point $x$, and killed at the boundary  
  of $C$. The distribution of this process at time $t$
   is therefore given by the
density,
   for $y\in C$,
   \begin{eqnarray*}
   p^0_t(x,y)\exp(\langle \xi,y-x\rangle-\frac{\Vert \xi\Vert^2}{2}t)
   &=&\sum_{w\in W}\varepsilon(w)
   p_t(x,w y)\exp(\langle \xi,y-x\rangle-\frac{\Vert \xi\Vert^2}{2}t)\\
   &=&\sum_{w\in W}\varepsilon(w)
   p_t(0,y-w x)\exp(\langle \xi,y-x\rangle-\frac{\Vert \xi\Vert^2}{2}t)
   \end{eqnarray*}
   where we have used the invariance of $p_t$ under the Weyl group.
   We now integrate this density over $C$,
    in order to get the probability that the
   exit time from $C$ is larger than $t$. Denoting by $T_C$ this exit
time, one has
   $$P(T_C>t)=\sum_{w\in W}\varepsilon(w)\int_C p_t(0,y-w x)
   \exp(\langle \xi,y-x\rangle-\frac{\Vert \xi\Vert^2}{2}t)\,dy$$
   Since the drift $\xi$ is in the chamber, for large $t$ one has
   $$\int_{V\setminus C}\exp(\langle \xi,y-x\rangle-\frac{\Vert
   \xi\Vert^2}{2}t)\,dy\to 0$$
    therefore
   $$\int_C p_t(0,y-w x)
   \exp(\langle \xi,y-x\rangle-\frac{\Vert
\xi\Vert^2}{2}t)\,dy\to_{t\to\infty}
   \exp(\langle \xi,w(x)-x\rangle)$$ and
   $$\lim_{t\to\infty}P(T_C>t)=P(T_C=\infty)=\sum_{w\in W}\varepsilon(w)
   \exp(\langle \xi,w(x)-x\rangle)$$
   We denote by $h_{\xi}(x)$ this function. It follows that,
conditionally on
   $\{T_C=\infty\}$, the Brownian motion with drift $\xi$, starting in  
$C$ and killed at the boundary of $C$, is a Markov process with
transition   probabilities
   $$q_{t,\xi}(x,y)
   =p^0_t(x,y)\frac{h_{\xi}(y)}{h_{\xi}(x)}\exp(\langle \xi,y-x\rangle
   -\frac{\Vert \xi\Vert^2}{2}t).$$
   Observe that $\frac{h_{\xi}(y)}{h_{\xi}(x)}\to\frac{h(y)}{h(x)}$ as
$\xi\to 0$.
   Standard arguments now show that as $x\to 0$ and $\xi\to 0$ the
distribution of
   this process approaches that of the Brownian motion in the Weyl  
chamber,
starting
   from 0.

   Finally we can rephrase this in the following way.
   \begin{lemma}\label{drift}
   The distribution of the Brownian motion with a drift $\xi\in C$,  
started
at 0 and
    conditioned to
   stay forever 
   in the cone $C-x$ (where $x\in C$) converges towards the
distribution of
  the Brownian motion in the Weyl chamber when $x,\xi\to 0$.
   \end{lemma}

   \subsection{Some further path transformations}

Let $w_0=s_{1}\ldots s_{q}$ be  a reduced decomposition and write
$\alpha_i=\alpha_{s_i}$.
Let $\eta:[0,T]\to V$ be a path with $\eta(0)=0$.   Recall that $\eta$
is dominant if
$\eta(t)\in\overline{C}$ for all $t\leq T$.  Set $\eta_q=\eta$ and,
for $j=1,\ldots ,q$,
$$\eta_{j-1}=\PP_{s_j}\ldots \PP_{s_q}\eta_q ,\qquad
x_j=-\inf_{0\leq t\leq T}\alpha_j^{\vee}( \eta_j(t)).$$
Then
$$\PP_{w_0}\eta(T)=\eta(T)+\sum_{j=1}^q x_j\alpha_j$$
and $\eta$ is dominant if, and only if, $x_j=0$ for all $j\le q$.
We now introduce some new path transformations and
give an alternative characterisation of dominant paths.

Let $w\in W$ be a  reflection, i.e.  $w$ is conjugate to some element in 
 $S$.  We choose   a non zero element $\alpha$
of
$V$ such that
$w\alpha=-\alpha$, then $w$ is the reflection $s_\alpha$ given by  
(\ref{salpha}) where $\alpha^\vee(v)=2(\alpha,v)/(\alpha,\alpha)$. As  
in
\cite{humphreys} we call $\alpha$ a positive root when $\alpha^\vee$ is  
positive on the Weyl chamber $C$, it is a simple  root when
$s_
\alpha
\in S$. Observe that one has
 $\PP_{\alpha}=\PP_{s_{\alpha}}$ for all positive 
roots (the left hand side is defined by Definition \ref{pitman-transform}, and
the second by Matsumoto's lemma, since $s_{\alpha}\in W$).

Let
$\beta$ be  a positive root, and
$s_{\beta}$ the associated  reflection.
For any positive root $\alpha$, one has $$s_{\beta}\, \pa\,
s_{\beta}=\mathcal P_{s_{\beta}(\alpha)}.$$
Consider the transformation
$\mathcal Q_{\beta}=\mathcal P_{\beta}\, s_{\beta}$.  One has
$$\mathcal Q_{\beta}\,\eta(t)=s_{\beta}\,\eta(t)+\sup_{0\leq s\leq
t}\beta^{\vee}(\eta(s))\,\beta .$$
Furthermore if $w_0=s_{1}\ldots s_{q}$ is a reduced decomposition  
$(s_i\in S)$, then
$$\mathcal Q_{w_0}:= \mathcal P_{w_0}\, w_0=\mathcal
Q_{\beta_{1}}\,
\ldots \, \mathcal Q_{\beta_{q}}$$
where $\beta_1=\alpha_1$ and $\beta_j=s_1\ldots s_{j-1}\alpha_j$.

Now define transformations $\da=s_{\alpha}\ea=
\iota\, \qa \, \iota$, where $\iota=-\kappa$.
  One has
\begin{equation}\label{da}
\D_\alpha\eta(t) = \eta(t)+\inf_{T\ge u\ge
t}\alpha^\vee(\eta(u)-\eta(t))\alpha
- \inf_{T\ge u\ge 0}\alpha^\vee(\eta(u))\alpha .
\end{equation}
Set $$\D_{w_0}= \D_{\beta_1} \cdots  \D_{\beta_q} =w_0\,\ewo=
\iota\,\q_{w_0}\,\iota $$
and note that $\D_{w_0} = \iota\,\PP_{w_0}\,(-\kappa)\, w_0$.

For a path $\eta$, 
set $\rho_q=\eta$ and , for $j\le q$,
$$\rho_{j-1}=\D_{\beta_j}\ldots\D_{\beta_q}\,\rho_q , \qquad
y_j=-\inf_{T\ge u\ge 0}\beta_j^{\vee}( \rho_j(u)).$$
\begin{lemma}\label{lem-d}
For all paths $\eta$ one has
\begin{equation}\label{eq-y}
\PP_{w_0}\eta(T)=\eta(T)+\sum_{j=1}^q y_j\beta_j .
\end{equation}
In particular, $\eta$ is dominant if, and only if, $y_j=0$ for all
$j\le q$.
\end{lemma}
{ Proof.}
By construction,
$$\mathcal D_{w_0}\eta (T)=\eta(T)+\sum_{j=1}^q y_j\beta_j .$$
Since $\mathcal D_{w_0}\eta(T)
=\PP_{w_0}\eta(T)$ by  Proposition \ref{e-duality},
this implies (\ref{eq-y}).  The path $\eta$ is dominant if, and only if,
$\PP_{w_0}\eta(T)=\eta(T)$.  By (\ref{eq-y}), this holds if, and only
if,
$\sum_j y_j\beta_j=0$ and, since the $y_j$ and $\beta_j$ are all
positive,
this is equivalent to the statement that $y_j=0$ for all $j\le q$.
\hfill $\diamondsuit$

\subsection{The representation theorem, first proof}

The definitions of transformations $\pa$, $\mathcal P_{w_0}$, $\qa$,
$\mathcal Q_{w_0}$ extend naturally to paths
$\pi$ defined on $\mathbb R^+$.
In this section we will prove  that, if $X$ is a
Brownian motion in $V$ (started from the origin),
then $\mathcal Q_{w_0}X$ is a Brownian motion in the fundamental
chamber $\overline C$.
Since $w_0$ leaves the distribution of Brownian motion invariant, this
implies that
$\mathcal P_{w_0}X$ is a Brownian motion in $\overline C$.

To prove this, we first extend the definition of the $\db$.  Let
$\beta$ be a positive root.
For paths $\pi:[0,+\infty)\to V$ with $\pi(0)=0$ and
$\alpha^\vee(\pi(t))\to+\infty$
as $t\to+\infty$ for all simple roots $\alpha$, define
\begin{equation}\label{D-def}
\mathcal D_{\beta}\pi(t)=\pi(t)+\inf_{ s\geq
t}\beta^{\vee}(\pi(s)-\pi(t))\beta
-\inf_{s\geq 0}\beta^{\vee}(\pi(s))\beta .
\end{equation}
Now set $\mathcal D_{w_0}= \mathcal D_{\beta_1} \cdots 
\mathcal D_{\beta_q}$
as before.   Since  $\mathcal D_{w_0}$ does not depend on the chosen
reduced
decomposition of $w_0$ we can also write
$\mathcal D_{w_0}= \mathcal D_{\beta_q} \cdots  \mathcal
D_{\beta_1}$.

\begin{lemma} \label{guinness}
If $\pi$ is a dominant path, one has
$\mathcal Q_{w_0}\, \mathcal D_{w_0}\,\pi =\pi$.
 \end{lemma}
{ Proof.}  It is easy to see that for any positive root $\beta$ and
path
$\xi:[0,\infty)\to V$ with $\xi(0)=0$ and
$\inf_{t\ge 0}\beta^\vee(\xi(t))=0$ we have $\mathcal Q_{\beta}
\mathcal D_{\beta}\xi =\xi$.
Let $\eta_0=\pi$ and $$\eta_j=\mathcal D_{\beta_j}\ldots \mathcal
D_{\beta_1}\pi,\qquad
v_j(t):=-\inf_{u\geq t}\beta^{\vee}_j(\eta_{j-1}(u)-\eta_{j-1}(t)).$$
Since $\pi$ is dominant
we have, by lemma \ref{lem-d} (with $T\to\infty$) that $v_j(0)=0$ for
each $j\le q$ and hence
$$\mathcal Q_{w_0}\, \mathcal D_{w_0}\pi =
\mathcal Q_{\beta_{1}}\ldots \mathcal Q_{\beta_{q}}\,
\mathcal D_{\beta_{q}}\ldots \mathcal D_{\beta_{1}}\,
\pi = \pi$$
as required.  \hfill $\diamondsuit$

\begin{lemma}\label{indD} If $X$ is a Brownian motion with drift in $C$, then
$\mathcal D_{w_0}X$
has the same distribution as $X$ and, moreover, is independent of the
collection of random
variables $\{ \inf_{t\ge 0} \alpha^\vee(X(t)),\ \alpha\mbox{
simple root}\}$.\end{lemma}
{ Proof.}
To prove this, we first need to extend the definitions of $\mathcal
D_\beta$ and
$\mathcal Q_\beta$ to paths $\pi$ defined on $\mathbb R$ with
$\pi(0)=0$ and
$\alpha^\vee(\pi(t))\to\pm\infty$ as $t\to \pm\infty$ for all simple
$\alpha$.
For $t\in\mathbb R$, set
$$\mathcal Q_{\beta}\pi(t)=s_{\beta}\,\pi(t)+\sup_{s\leq
t}\beta^{\vee}(\pi(s))\,\beta
-\sup_{s\leq 0}\beta^{\vee}(\pi(s))\,\beta $$
and define $\mathcal D_\beta \pi$ by (\ref{D-def}) allowing
$t\in\mathbb R$.
Then, if $\iota$ denotes the involution $$\iota\,\pi(t)=-\pi(-t)$$
one has $\mathcal D_{\beta}=\iota\, \mathcal Q_{\beta}\, \iota$
and $\mathcal D_{w_0}:={\mathcal D}_{\beta_q}\,\cdots
\,
\mathcal D_{\beta_1} = \iota\, \q_{w_0}\,\iota$ as before.
Note that $\D_{w_0}$ does not depend on the particular reduced
decomposition of $w_0$,
and also that $\mathcal D_\beta (\pi(t),t\ge 0)=(\mathcal D_\beta
\pi(t),t\ge 0)$
and $\mathcal D_{w_0} (\pi(t),t\ge 0)=(\mathcal D_{w_0} \pi(t),t\ge 0)$.
We will use the following auxillary lemma.
\begin{lemma}\label{V}
Let $\pi: \mathbb R \to V$ with $\pi(0)=0$, and
$\alpha(\pi(t))\to\pm\infty$ as $t\to\pm \infty$ for all simple roots
$\alpha$.
Then, for all $t\in \mathbb R$,
$$-\inf_{u\geq t}\beta^{\vee}(\pi(u)-\pi(t))=-\inf_{s\leq
t}\beta^{\vee}(\db\,\pi(u)-\db\,\pi(t)).$$
\end{lemma}
{ Proof.}
This can be checked directly, or deduced from (\ref{eta-pi}).
\hfill$\diamondsuit$
\medskip

Introduce a Brownian motion $Y$ indexed by $\mathbb R$ such that
$X=(Y(t),t\ge 0)$ and $(\iota Y(t),t\ge 0)$ is an independent copy of
$X$.
For any positive root $\beta$, the distribution of $\db Y$ is the same
as that of $Y$.
This is a one-dimensional statement which can be checked directly,
or can be seen as a consequence of the classical output
theorem on the $M/M/1$ queue (see, for example,
\cite{oy1}).
In particular, the distribution of $\db X$ is the same as that of $X$.
It follows that $\mathcal D_{w_0}Y$ has the same distribution as $Y$,
and
$\mathcal D_{w_0}X$ has the same distribution as $X$.
Let $Y_0=Y$ and
   $$Y_j=\mathcal D_{\beta_j}\ldots \mathcal D_{\beta_1}Y,\qquad
V_j(t):=-\inf_{u\geq t}\beta^{\vee}_j(Y_{j-1}(u)-Y_{j-1}(t)).$$
Note that $Y_q=\mathcal D_{w_0}Y$ and recall that, for $t\ge 0$,
$\mathcal D_{w_0}Y(t)=\mathcal D_{w_0}X(t)$.
By  Lemma \ref{V} one has
$$V_q(t)=-\inf_{s\leq t}\beta^{\vee}_j(Y_{q}(s)-Y_{q}(t))$$
$$Y_{q-1}(t)=Y_q(t)+(V_q(t)-V_q(0))\beta_q$$
and by induction on $k$,
$$V_{q-k}(t)=-\inf_{s\leq t}\beta^{\vee}_j(Y_{q-k}(s)-Y_{q-k}(t))$$
$$Y_{q-k-1}(t)=Y_{q-k}(t)+(V_{q-k}(t)-V_{q-k}(0))\beta_{q-k}$$
It follows that the $(V_j(t),t\leq 0)$ are measurable with respect to the
$\sigma$-field
generated by $(\mathcal D_{w_0}Y(s),s\leq 0)$.  In particular, the random
variable
$V_1(0)=\inf_{t>0}\beta_1^\vee(X(t))$ is measurable with respect to the
$\sigma$-field
generated by $(\mathcal D_{w_0}Y(s),s\leq 0)$.  Now, for each $\alpha\in
S$, there
is a reduced decomposition of $w_0$ with $\beta_1=\alpha$, so we see
that
the random variables $\{ \inf_{t\ge 0} \alpha^\vee(X(t)),\ \alpha\mbox{
simple}\}$
are all measurable with respect to the $\sigma$-field generated by
$(\mathcal D_{w_0}Y(s),s\leq 0)$,
and therefore independent of $(\mathcal D_{w_0}Y(s),s\geq 0)$, as
required. $\hfill\diamondsuit$

\begin{theorem}\label{rep_thm}
Let $X$ be a Brownian motion in $V$. Then 
$\mathcal P_{w_0} X$ is a Brownian motion in $\overline C$.
\end{theorem}
{ Proof.}
Let $x,\xi\in C$ and let $X$ be a Brownian motion with drift $\xi$.
   The event `$X$ remains in the cone $C-x$ for all
times' can be expressed in terms of the variables
$\{ \inf_{t\ge 0} \alpha^\vee(X(t)),\ \alpha\mbox{
simple root}\}$
 therefore, by 
 lemma \ref{indD}, it 
is independent of $(\mathcal D_{w_0}X(t),t\geq 0)$.  Thus, if $R$ has the
same
distribution as that of $X$ conditioned on this event, then $\mathcal
D_{w_0}R$
has the same distribution as $X$.  Now we can let $x,\xi\to 0$ so that
$X$ is a
Brownian motion with no drift and $R$ is a Brownian motion in
$\overline C$;
by continuity, $\mathcal D_{w_0}R$ has the same distribution as $X$.
Now, by lemma \ref{guinness} $\mathcal Q_{w_0} \mathcal D_{w_0}R = R$
almost surely.  It follows that $\mathcal Q_{w_0} X$, and hence
$\mathcal P_{w_0} X$,
is a Brownian motion in $\overline C$, as required.
\hfill $\diamondsuit$

\subsection{ Random walks and Markov chains on the weight  
lattice}\label{minus}
We will now present the second proof of the Brownian motion property.
We assume  that $W$ is the Weyl group of the semisimple  
Lie algebra $\mathfrak g$ as in sections \ref{sec_semisimple},
\ref{BMWC}, and $V=\mathfrak a^*$.  As in section \ref{BMWC}, let
$T$ be a maximal torus of the compact group $K$, the simply connected  
compact group with Lie algebra
  $\mathfrak g_{\mathbb R}$, a compact form of $\mathfrak g$.  Let
$\omega\in P_+$ be a nonzero dominant weight and let $\chi_\omega$ be  
the  character
of the  associated highest weight module.
As a function on  $T$ this is the Fourier transform of the
positive measure $R_{\omega}$ on $P$,  which puts a weight
  $m_{\mu}^{\omega}$ on a weight
    $\mu$ where $m_{\mu}^{\omega}$ is
the multiplicity of  $\mu$ in the module with highest weight $\omega$.  
In other
words
$$\chi_{\omega}=\sum_{\mu\in P_+}m^{\omega}_{\mu}e(\mu)$$
where $e(\mu)(\theta)=e^{2i\pi\langle\mu,\theta\rangle}$ is the  
character on $T$.
We can
divide this measure $R_\omega$ by  $\text{dim}\,\omega$ to get a  
probability measure
$$\nu_\omega=\frac{1}{\dim \omega}R_{\omega}.$$
Consider  the random walk ($X_n,n\geq 0$),
on the weight lattice, whose increments are
distributed according to this probability measure, started at zero.
Thus the transition probabilities of this random walk are given by
$$p_{\omega}(\mu,\lambda)=\frac{m^{\omega}_{\lambda- 
\mu}}{\text{dim}\,{\omega}}.$$
  Donsker's
theorem and invariance of $m^\omega$ under the Weyl group implies
\begin{theorem}\label{donsker} The stochastic process
$\frac{X_{[Nt]}}{\sqrt{N}}$ converges, as $N\to\infty$, to a Brownian  
motion on
$\mathfrak a^*$ with  correlation invariant under $W$.
\end{theorem}

 Let us define  a probability  
transition
  function $q_{\omega}$ on $P_+$ by the formula
   
$$\frac{\chi_{\mu}}{\text{dim}\,\mu}\frac{\chi_{\omega}}{\text{dim}\,\omega}
  =\sum_{\lambda\in
  P_+}q_\omega(\mu,\lambda)\frac{\chi_{\lambda}}{\text{dim}\,\lambda}.$$
Thus $q_\omega(\mu,\lambda)$ is equal to
  $\frac{M_{\omega,\mu}^{\lambda}\text{dim}\,\lambda}
  {\text{dim}\,\omega\,\text{dim}\,\mu}$ where  
$M_{\omega,\mu}^{\lambda}$ is the
  multiplicity of the module with  highest weight
   $\lambda$ in the decomposition of the
  tensor product of the modules with highest weights $\omega$ and $  
\mu$, see, e.g. \cite{Eymard},\cite{bi2}.
  \begin{lemma}\label{qomega}
  One has
   $$
   q_{\omega}(\mu,\lambda)=\frac{\text{\rm dim}\,\lambda}{\text{\rm  
dim}\,
  \mu}
  \sum_{w\in W}\varepsilon (w)p_{\omega}(\mu+\rho,w(\lambda+\rho)).
  $$
  \end{lemma}
    { Proof.}  Let $dk$ be the normalized Haar measure
  on $K$. By the
  orthogonality relations for characters,   one has
  $$M_{\omega,\mu}^{\lambda}=
  \int_{K}\overline{\chi_{\lambda}}(k)\chi_{\mu}(k)
  \chi_{\omega}(k)dk$$
  therefore
   
$$q_{\omega}(\mu,\lambda)=\frac{M_{\omega,\mu}^{\lambda}\text{dim}\,\lambda}
   {\text{dim}\,\omega\,\text{dim}\,\mu}=\frac{\text{dim}\,\lambda}{\text{dim}\,
  \mu\,\text{dim}\,\omega}
  \int_{K}\overline{\chi_{\lambda}}(k)\chi_{\mu}(k)
  \chi_{\omega}(k)dk$$
  Now we can use the Weyl integration formula as well as Weyl's  
character formula
  to rewrite the formula as an integral over $T$, the maximal torus of  
$K$.
  Thus
  $$q_{\omega}(\mu,\lambda)=\frac{|W|\,\text{dim}\,\lambda}{\text{dim}\,
  \mu\,\text{dim}\,\omega}
  \int_{T}\sum_{w_1,w_2\in W}\varepsilon(w_1)\varepsilon(w_2)
  \overline{e(w_1(\lambda+\rho)}(\theta)e(w_2(\mu+\rho))(\theta)
  \chi_{\omega}(\theta)d\theta$$
  where $e(\gamma)(\theta)=e^{2i\pi\langle\gamma,\theta\rangle}$ and  
$\rho$ is
  half the sum of positive weights.
  Now using the invariance of $\chi_{\omega}$
  under the Weyl group we can rewrite this as
   \begin{eqnarray*}
   q_{\omega}(\mu,\lambda)&=&\frac{\text{dim}\,\lambda}{\text{dim}\,
  \mu\,\text{dim}\,\omega}
  \int_{T}\sum_{w\in W}\varepsilon(w)
  \overline{e(\lambda+\rho)}(\theta)e(w(\mu+\rho))(\theta)
  \chi_{\omega}(\theta)d\theta
  \\
  &=&\frac{\text{dim}\,\lambda}{\text{dim}\,
  \mu}
  \sum_{w\in W}p_{\omega}(\mu+\rho,w(\lambda+\rho)).
  \end{eqnarray*}
   \hfill $\diamondsuit$

    \medskip
    From (\ref{q-p}), Theorem \ref{donsker}, Lemma \ref{qomega},
     and standard arguments, we deduce

    \begin{prop}\label{Y-bmwc}
  Let $Y$ be a Markov chain on $P_+$ started at 0, with transition  
probabilities
   $q_\omega(\mu,\lambda)$, then $\frac{Y([Nt])}{\sqrt{N}}$ converges in
   distribution, as $N\to\infty$ to a Brownian motion in the Weyl chamber
   $\overline C$.
  \end{prop}
  \subsection{Pitman operators and the Markov chain on the weight  
lattice}
We  choose a nonzero dominant
weight $\omega$, and a dominant path $\pi^{\omega}$ defined on $[0,1]$  
with $\pi^{\omega}(1)=\omega$.
  Let $B\pi^{\omega}$ be
the set of paths in the Littelmann module generated by $\pi^{\omega}$.
  We now construct a
stochastic process with values in $P$. Choose independent  random
paths $(\eta_n\in B\pi^{\omega}, n=1,2,\ldots)$,
each with uniform distribution on $B\pi^{\omega}$,
  and define the stochastic process $Z$
as the random path  obtained by the usual concatenations  
$\eta_1*\eta_2*\cdots$ of the
  $\eta_i;i=1,2,\ldots$. In
other words, one has
$Z(t)=\eta_1(1)+\eta_2(1)+\ldots+\eta_{n-1}(1)+\eta_n(t-n)$ if
$t\in [n,n+1]$.
Beware
 that this concatenation does not coincide with Littelmann's definition,
recalled in section (\ref{LR}), since we do not rescale the time. 
Littelmann's theory then implies that  $\eta_n(1)$ is a
random weight in $P$ with distribution $
\nu_{\omega}$, and
  ($Z(n),n=0,1,\ldots$) is the
random walk in $\mathfrak a^*$ with this distribution of increments.
\begin{theorem}\label{rwwc}
The stochastic process $(\PP_{w_0}Z(n),n=0,1,\ldots)$ is a Markov chain  
on $P_+$,
with probability transitions $q_\omega$.
\end{theorem}
{ Proof.} First note that the set of
paths of the form  $\eta_1*\eta_2*\ldots *\eta_n$ where $\eta_i\in  
B\pi^{\omega}$
is stable under Littelmann
operators, by \cite{littel}, therefore by (\ref{pit/lit}) it is also
    stable under Pitman transformations.
Consider a dominant  path of the form
  $\gamma_1*\gamma_2*\ldots *\gamma_n$, with all $\gamma_i\in  
B\pi^{\omega}$.
  We shall compute the conditional
  probability distribution of
  $\PP_{w_0}Z(n+1)$ knowing that $\PP_{w_0}Z(t)
  =\gamma_1*\gamma_2*\ldots *\gamma_n(t)$
  for $t\leq n$.
   Let $\mu=\gamma_1*\gamma_2*\ldots *\gamma_n(1)$.
   By Corollary \ref{pit-littel} the set of all paths of the form
  $\eta_1*\eta_2*\ldots *\eta_n$ such that  
$\PP_{w_0}(\eta_1*\eta_2*\ldots *\eta_n)
  =\gamma_1*\gamma_2*\ldots *\gamma_n$ coincides with the Littelmann
  module $B(\gamma_1*\gamma_2*\ldots *\gamma_n)$. Now consider a path
   $\eta_{n+1}\in B\pi$ and the concatenation $\eta_1*\eta_2*\ldots
   *\eta_n*\eta_{n+1}$,
   then $\PP_{w_0}(\eta_1*\eta_2*\ldots *\eta_n*\eta_{n+1})$ will be the  
dominant
   path in the Littelmann module generated by
   $\eta_1*\eta_2*\ldots *\eta_n*\eta_{n+1}$. By Littelmann's version of  
the
   Littlewood-Richardson rule (section 10 in \cite{littel}), the number  
of pairs of paths
   $(\eta_1*\eta_2*\ldots *\eta_n,\eta_{n+1})$  such that
   $\PP_{w_0}(\eta_1*\eta_2*\ldots *\eta_n)
  =\gamma_1*\gamma_2*\ldots *\gamma_n$ and
  $\PP_{w_0}(\eta_1*\eta_2*\ldots *\eta_n*\eta_{n+1})(1)=\lambda$ is  
equal to the
  dimension of the isotypic component of type $\lambda$ in the module  
which is
  the tensor product
  of the highest weight modules $\mu$ and $\omega$, in particular this  
depends
  only on $\mu$, and is equal to
  $M_{\omega,\mu}^{\lambda}\text{dim}\,\lambda$. Since the total number  
of pairs
  $(\eta_1*\eta_2*\ldots *\eta_n,\eta_{n+1})$ with
  $\PP_{w_0}(\eta_1*\eta_2*\ldots *\eta_n)
  =\gamma_1*\gamma_2*\ldots *\gamma_n$ is $\text{dim}\,  
\mu\,\text{dim}\,\omega$, we see
  that the conditional probability we seek is
  $\frac{M_{\omega,\mu}^{\lambda}\text{dim}\,\lambda}
  {\text{dim}\,\omega\,\text{dim}\,\mu}=q_\omega(\mu,\lambda)$.
  This proves the claim.
\hfill $\diamondsuit$

\subsection{Second proof of the representation theorem for Weyl groups}
Putting together Proposition \ref{Y-bmwc} and Theorem \ref{rwwc} we get  
another proof of
Theorem \ref{rep_thm}. Indeed, by Donsker's theorem, the process   
$\frac{Z([Nt])}{\sqrt{N}}$
  gives as limit the Brownian motion in $\mathfrak a^*$.
The process $(\mathcal P_{w_0}Z(n),n\geq 0)$ is distributed as the  
Markov process
   of
  Proposition \ref{Y-bmwc}, by Theorem \ref{rwwc}.
   Applying the scaling of Proposition \ref{Y-bmwc} to the stochastic
    process
   $(\mathcal P_{w_0}Z(t),t\geq 0)$ yields for limit process  the
   Brownian motion on the Weyl chamber. Since $\mathcal P_{w_0}$ is a  
continuous
   map, which commutes with
   scaling we get the proof of Theorem \ref{rep_thm}, when $W$ is the  
Weyl group of a complex semisimple Lie algebra.\hfill $\diamondsuit$

\subsection{A remark on the Duistermaat-Heckman measure} 
The distribution of the path $t\in  
[0,n]\mapsto Z(t)$ is uniform on the set $$B(\pi^{\omega})^{*n}=
\{\eta_1*\eta_2*\cdots*\eta_n; \eta_i\in B\pi^\omega\}.$$

Therefore, for any path $\eta \in B(\pi^{\omega})^{*n}$, the  
distribution of $(Z(s))_{0 \leq s \leq n}$ conditionally on  
$\{\PP_{w_0}Z(s)=
\eta(s), 0 \leq s \leq n\},$
is uniform on the set $\{\gamma \in B(\pi^{\omega})^{*n};  
\PP_{w_0}\gamma=\eta\}$. It thus follows from Littelmann theory  
\cite{littel} that
the conditional distribution of the terminal value $Z_n$ is the  
probability measure $\nu_\eta$. It has been proved by Heckman  
\cite{hec} (see
also \cite{gui}, \cite{Dui}) that if $\gamma_\varepsilon \to \infty$ in  
$\mathfrak a^*_+$ and $\varepsilon \gamma_\varepsilon \to v $ then
$D_\varepsilon\nu_{\gamma_\varepsilon}$ converges to the so called  
Duistermaat-Heckman measure associated to $v$, i.e. the projection of  
the
normalized measure on the coadjoint orbit of $K$ through $v$, by the  
orthogonal projection on $\mathfrak a^*$. 
This follows from Kirillov's character formula for $K$.
From the preceding
section we deduce that if
$X$ is the Brownian motion on $\mathfrak a^*$, then the  law of $X(T)$  
conditionally on $\PP_{w_0}X=\gamma$ on $[0,T]$ is the 
Duistermaat-Heckman measure
associated with $\gamma(T)$.

    \section{Appendix.
    Proof of Proposition \ref{pit} ({\it iv})}\label{prfpit}
    Let $\eta$ be a path. Defining $\pi=\pa\eta$,  $x=-\inf_{T\geq t\geq
    0}\alpha^{\vee}(\eta(t))$, and  $t_0=\sup\{t|\alpha^{\vee}(\eta(t))=-x\}$, we
shall
    check that equation \ref{eta-pi} is valid.

    If $t\geq t_0$ then one has
    $\inf_{0\leq s\leq t}
\alpha^{\vee}( \eta(s) )=-x $ therefore
    \begin{eqnarray*}\alpha^{\vee}(\pi(t))&=&\alpha^{\vee}(\eta(t))+2x\\
    &=&x+(\alpha^{\vee}(\eta(t))+x)
    \\&\geq& x\end{eqnarray*}
    for all $t\geq t_0$. It follows that
    $\inf\left(x,\inf_{T\geq s\geq t}
\alpha^{\vee}( \pi(s) )\right)=x$ for $t\geq t_0$.
    Formula \ref{eta-pi} follows for $t\geq t_0$.

    If $t<t_0$, let $u=\inf\{s\geq t|\alpha^{\vee}(\eta(s))=\inf _{0\leq
v\leq
    t}\alpha^{\vee}(\eta(v))\}$. Then
    $t\leq u\leq t_0$.
    One has \begin{eqnarray*}
    \alpha^{\vee}(\pi(u))&=&\alpha^{\vee}(\eta(u))-2\inf_{0\leq v\leq
    u}\alpha^{\vee}(\eta(v))\\&=&-\alpha^{\vee}(\eta(u))
    \end{eqnarray*}
    which implies that $\inf_{T\geq v\geq t}\alpha^{\vee}(\pi(v))\leq
    -\inf _{0\leq v\leq
    t}\alpha^{\vee}(\eta(v))\}\leq x$.
    On the other hand, for  $v\geq t$ one has
    \begin{eqnarray*}
    \alpha^{\vee}(\pi(v))&=&\alpha^{\vee}(\eta(v))-2\inf _{0\leq s\leq
    v}\alpha^{\vee}(\eta(s))\\
    &\geq&(\alpha^{\vee}(\eta(v))-\inf _{0\leq s\leq
    v}\alpha^{\vee}(\eta(s)))-\inf _{0\leq s\leq
    t}\alpha^{\vee}(\eta(s))
    \\&\geq&-\inf _{0\leq s\leq
    t}\alpha^{\vee}(\eta(s))
    \end{eqnarray*}
    therefore
    $\inf_{T\geq v\geq t}\alpha^{\vee}(\pi(v))=-\inf _{0\leq s\leq
    t}\alpha^{\vee}(\eta(s))$ and Formula \ref{eta-pi} for $t<t_0$
follows.
    The existence and uniqueness in Proposition \ref{pit} follows.
  \hfill$\diamondsuit$


\end{document}